\documentclass[final]{siamltex}
\usepackage{latexsym,url,amscd,amsmath,amssymb}
\usepackage{graphics}
\usepackage{graphicx}
\usepackage{epstopdf}
\usepackage[usenames]{color}
\usepackage{hyperref}

\newtheorem{assumption}{Assumption}

\newtheorem{algorithm}{Algorithm}

\newcommand{\mat}[1]{\left[ \begin{array}{#1} }
\newcommand{\rix}{\end{array} \right]}

\newcommand{\vw}{\mathbf{w}}

\newcommand{\R}{\mathbb{R}}

\usepackage[usenames]{color}

\begin{document}

\title{A Fast Algorithm for Total Variation Image Reconstruction \\
from Random Projections}

\author{Yunhai Xiao\footnotemark[1] \and Junfeng Yang\footnotemark[2]}

\renewcommand{\thefootnote}{\fnsymbol{footnote}}
\footnotetext[1]{Institute
of Applied Mathematics, College of Mathematics and Information
Science, Henan University, Kaifeng, 475004, P. R. China. Current
Position: Department of Mathematics, Nanjing University, Nanjing,
210093, P. R. China (yhxiaomath@gmail.com). This author's work was
supported by the Chinese Post-doctoral Foundation grant 20090461094.}
\footnotetext[2]{Department of Mathematics, Nanjing
University, Nanjing, 210093, P. R. China (send correspondence to this author at:
jfyang@nju.edu.cn). This author's work was supported in part by the Natural
Science Foundation of China grant NSFC-10971095 and the Natural Science
Foundation of Jiangsu Province BK2008255.}

\renewcommand{\thefootnote}{\arabic{footnote}}

\date{\today}
\maketitle

\begin{abstract}
Total variation (TV) regularization is popular in image restoration
and reconstruction due to its ability to preserve image edges. To
date, most research activities on TV models concentrate on image
restoration from blurry and noisy observations, while discussions on
image reconstruction from random projections are relatively fewer.
In this paper, we propose, analyze, and test a fast alternating
minimization algorithm for image reconstruction from random
projections via solving a TV regularized least-squares problem.
The per-iteration cost of the proposed algorithm involves a linear time shrinkage operation, two matrix-vector multiplications and two fast Fourier transforms. Convergence, certain finite convergence and $q$-linear convergence results are established, which indicate that the asymptotic convergence speed of the proposed algorithm depends on the spectral radii of certain submatrix. Moreover, to speed up convergence and enhance robustness, we suggest an accelerated scheme based on an inexact alternating direction method. We present experimental results to compare with an existing algorithm, which indicate that the proposed algorithm is stable, efficient and competitive with TwIST \cite{TWIST} --- a state-of-the art algorithm for solving TV regularization problems.
\end{abstract}

\begin{keywords}
Total variation, image restoration, image reconstruction, compressive
sensing, alternating direction method
\end{keywords}

\begin{AMS}
68U10, 65J22, 65K10, 65T50, 90C25
\end{AMS}

\pagestyle{myheadings} \thispagestyle{plain} \markboth{Y.-H. Xiao
and J.-F. Yang}{A Fast Algorithm for TV Image Reconstruction
from Random Projections}

\section{Introduction}
Image restoration and reconstruction play  important roles in medical and
astronomical imaging,  image and video coding, file restoration, and many
other applications. Let $\bar{u}\in\R^{n^2}$ be an original $n\times n$ image,
$A\in \R^{m\times n^2}$ be a linear operator, and $f\in\R^m$ be an observation which satisfies the
relationship
\begin{equation}\label{obser}
f=\mathbb{N}(A\bar{u}) \in\R^m,
\end{equation}
where $\mathbb{N}(\cdot)$ represents a noise contamination or corruption procedure.
Given $A$, image restoration and reconstruction extract $\bar{u}$ from $f$,
which is either under-determined ($m<n^2$) or ill-possed (e.g.,  deconvolution/deblurring),
making classical least-squares approximation alone not suitable.
To stabilize recovery, regularization technique is frequently used, giving a general
reconstruction model of the form
\begin{equation}\label{P0-reg}
\min_{u}\Phi_{\text{reg}}(u)+ \mu\Phi_{\text{fid}}(Au-f),
\end{equation}
where $\Phi_{\text{reg}}(u)$ promotes solution regularity such as smoothness and sparseness, $\Phi_{\text{fid}}(Au-f)$ fits the observed data by penalizing
the difference between $Au$ and $f$, and $\mu>0$ balances the two terms for minimization. The choice of $\Phi_{\text{fid}}(\cdot)$ depends on different noise, e.g.,
the squared $\ell_2$ penalty is usually used for Gaussian additive noise, while the $\ell_1$ penalty is more appropriate for certain non-Gaussian noise, e.g.,
salt-and-pepper noise. Throughout this paper, we assume that $\mathbb{N}(\cdot)$ represents an additive Gaussian noise contamination and thus set $\Phi_{\text{fid}}(\cdot)=\|\cdot\|_2^2$.
Among other  regularization, total variation (TV) has been popular ever since its introduction by Rubin, Osher
and Fatemi \cite{TV}. The remarkable property of TV is to preserve edges due to its linear penalty on differences between adjacent pixels.
The most widely studied TV model for image deconvolution (in which case $A\in\R^{n^2\times n^2}$ is a convolution matrix) is
\begin{equation}\label{tv}
\min_{u}\sum\nolimits_{i=1}^{n^2}\|D_iu\|_2+\frac{\mu}{2}\|Au-f\|^2_2,
\end{equation}
where $D_i\in\R^{2\times n^2}$ denotes the local finite difference operator
(with ceratin boundary conditions) at pixel $i$, and $\sum_i\|D_iu\|_2$
is a discretization of the TV of  $u$. In this paper, we propose, analyze and test
a fast alternating minimization algorithm for solving \eqref{tv}, in which $A$ is a compressive sensing encoding matrix
($m < n^2$) and does not have    structures.

In the following of this section, we review briefly compressive sensing ideas and algorithms, which provide theoretical
guarantee for image reconstruction via solving  \eqref{tv}, examine some existing algorithms for relevant TV problems,
and describe  the contributions and organization of this paper. Throughout this paper, we refer to \eqref{tv} as TV/L$^2$.

\subsection{Compressive sensing --- ideas and algorithms}

Compressive sensing (CS) is an emerging methodology in digital signal processing brought to the research forefront by Donoho \cite{DONOBO1}, Cand{\` e}s, Romberg and Tao \cite{CANDES1,CANDES2}, and has attracted intensive research activities in the past few years. In a nutshell, CS first encodes a sparse signal (possibly under certain sparsifying basis) through hardware devices into a relatively small number of linear projections and then reconstructs it from the limited measurements. Let $\bar{x}\in\R^n$ be the sparse signal that we wish to capture, i.e., the number of nonzeros in $\bar x$ is much less than its length $n$, and $b=A\bar{x}$ represent a set of $m$ (usually much smaller than $n$) linear projections of $\bar{x}$. Under certain
desirable conditions, it is shown that  with high probability  the basis pursuit problem
\begin{equation}\label{basisp}
\min_x\|x\|_1\quad \mbox{s.t.} \ Ax=b,
\end{equation}
yields the sparsest solution of the linear system $Ax=b$, see \cite{SHARE}.
More often than not, $b$ contains noise, in which case certain relaxation is desirable.
For white Gaussian noise, the most widely used models are the basis pursuit denoising problem
\begin{equation}\label{csone}
\min_{x}\|x\|_1+\frac{\mu}{2}\|Ax-b\|_2^2,
\end{equation}
where, roughly speaking, $\mu>0$ is inversely proportional to the noise level, and its variants.
Since most signals of interests are sparse or nearly sparse (called compressible) under certain basis,
the CS idea has extremely wide applications.  Recent results
show  that stable reconstruction can be obtained provided that $A$
possesses certain randomness. It has been clear from \cite{YZHANG1} that
for almost all random matrices the exact recoverability is approximatively
identical. Moreover, exact recoverability is attainable when $A$ contains randomly taken
rows from orthonormal matrices, e.g., partial Fourier which arises  from magnetic
resonance imaging \cite{MRI}.

In the application of CS, matrix $A$ is large and dense. Furthermore, in certain applications $A$ contains structures
that allow fast matrix-vector multiplication, e.g., $A$ is a partial Fourier matrix as in MRI.
These features make traditional powerful optimization approaches such as interior point methods not suitable.
In comparison, first-order algorithms that depend on merely matrix-vector multiplications are more desirable.
Therefore, in the last few years numerous algorithms have been proposed for recovering
sparse signals via solving certain $\ell_1$-norm regularized problems including \eqref{basisp}, \eqref{csone} and thier variants.
Several  well-known  approaches in this area include the gradient projection   method  \cite{FNW1},
the fixed-point continuation   method \cite{ZHANG2}, the spectral projected gradient   method  \cite{SPGL1}, and
the Bregman iterative method   \cite{OSHER2,WYIN2,WYIN3,CAI1,CAI2}. More recent algorithms can be found in \cite{AB,PLC,WEN,YANG5}.

\subsection{Some existing algorithms for TV/L$^2$}
The advantage of TV regularization compared with Tikhonov-like \cite{TIKHONOV} regularization
in recovering high quality image is not without a price. The nondifferentiability
of TV causes the main difficulty. In addition, problems arising  from signal and image reconstruction
are usually large scale and ill-possed, which further make TV models difficult to be solved efficiently.
Since the introduction of TV regularization, many algorithms have been proposed for solving \eqref{tv} and its variants.
In the pioneer work \cite{TV}, a time-marching scheme was used to solve a partial
differential equation system, which in optimization point of view is
equivalent to a constant step-length gradient descent method.
This time-marching scheme suffers slow convergence especially when the
iterate point approaches the solution set.  Another well-known method is
the linearized gradient method proposed in \cite{VCR96} for
denoising and in \cite{AF98} for deblurring, which solves the Euler-Lagrangian equation
via a fixed-point iteration. At each iteration of the linearized
gradient method, a linear system  needs to be solved, which makes the per-iteration cost
extremely expensive especially when the problem becomes
more ill-conditioned. To overcome the linear convergence of first-order methods, the
authors of \cite{VCR96} incorporated Newton method to solve
(\ref{tv}), which achieved superlinear convergence at the cost of solving a large linear system at each iteration.
Another important approach for TV problems is the iterative
shrinkage/thresholding (IST) method \cite{ME06,MEBM06,JLMN03}.
In \cite{TWIST}, Bioucas-Dias and Figueiredo introduced a two-step IST (TwIST)  algorithm, which
exhibits much faster convergence than the primary IST algorithm for
ill-conditioned problems. We note that IST-based algorithms require to solve a
TV denoising subproblem at each iteration  which requires its own iterations.

Despite the progress have been achieved, algorithms for solving
\eqref{tv} are still much slower than those for Tikhonov regularization problems.
Recently, a fast TV deconvolution (FTVd) method is proposed in \cite{YANG1}, which makes full use of
problem structures (both $A$ and finite difference operators have circulant structures under proper boundary conditions)
and thus converges very fast.  FTVd solves  a penalty approximation of \eqref{tv}, that is
\begin{equation}\label{atv}
\min_{u,\mathbf{w}}\sum\nolimits_{i=1}^{n^2}\left(\|\mathbf{w}_i\|_2+\frac{\beta}{2}\|\mathbf{w}_i-D_iu\|_2^2\right)
+\frac{\mu}{2}\|Au-f\|_2^2,
\end{equation}
where, for each $i$, $\mathbf{w}_i\in\R^2$ is an auxiliary variable and $\beta>0$ is a
penalty parameter. The advantage of considering \eqref{atv} is that it leads to fast and efficient
alternating minimizations for deconvolution problems. The numerical results given in
\cite{YANG1} indicates that FTVd is much faster than the lagged diffusivity method in \cite{VCR96}, which is known  to be efficient previously.
For more details on the FTVd algorithm and its performance, see \cite{YANG1}. Given the practical efficiency of FTVd,
this split and penalty idea has been  extended to multichannel image restoration
in \cite{YANG2}, impulsive noise elimination in \cite{YANG3} and medical reconstruction from partial
Fourier coefficients in \cite{YANG4}. More algorithms for TV/L$^2$ problem can be found in
\cite{TFCHAN1,TFCHAN2,NG3,NG1,PCF,NG4} and references therein.

\subsection{Contributions}

The purpose of this paper is to develop a fast algorithm for solving  \eqref{tv}, where $A$ is
a general linear operator. Specifically, we are interested
in compressive sensing encoding matrices in which case $A$ contains smaller or even much smaller
number of rows than columns. As is stated above, problem \eqref{atv} admits fast alternating minimization when $A$ is a convolution matrix.
As a matter of fact, the minimization of \eqref{atv} with respect to $\vw_i$, $i=1,2,\ldots,n^2$,
reduces to $n^2$ two-dimensional problems (no matter what $A$ is), which can be solved
easily and exactly in linear time. However, different from deconvolution problems, $A$ does
not have structures in our stated case. Consequently, the solution of $u$-subproblems can  not utilize any fast transforms.

In this paper, we first introduce a fast alternating minimization
scheme for solving \eqref{atv}, which recurs to linearization and proximal
techniques when solving the $u$-subproblems. Under quite reasonable
technical assumptions,  we show that the proposed algorithm
converges globally to a solution of \eqref{atv}. Moreover, we establish
$q$-linear convergence results which indicate  that the $q$-linear
factor  depends on the spectral radius of certain submatrix. Clearly,
the solution of \eqref{atv} well approximates that of (\ref{tv})
only when  $\beta$  is sufficiently large, which causes numerical
difficulties in computation. To overcome this drawback, we introduce an
inexact alternating direction method, which accelerates the
convergence of the alternating minimization approach and converges
to a solution of \eqref{tv} without driving $\beta$ to infinity.
Since the proposed algorithms solve \eqref{atv} and \eqref{tv} with a CS encoding matrix,
we name the resulting algorithms FTVCS.  We present experimental
results and compare with TwIST \cite{TWIST}. The comparison results
indicate that FTVCS is fast and efficient and performs comparable
with the state-of-the art algorithm TwIST.

\subsection{Notation and organization}
Now, we define our notation.
For scalars $\alpha_i$, vectors $v_i$, and matrices
$M_i$ of appropriate sizes, $i = 1, 2$, we let $\alpha = (\alpha_1;\alpha_2)\triangleq(\alpha_1,\alpha_2)^\top$,
$v = (v_1;v_2)\triangleq(v^\top_1, v^\top_2)^\top$, and $M = (M_1;M_2)\triangleq(M^\top_1,M^\top_2)^\top$.
Let $D^{(1)}$ and $D^{(2)}$ be the two first-order finite difference matrices in horizontal and
vertical directions, respectively. As is used before, $D_i\in \R^{2\times n^2}$ is a two-row matrix formed by
stacking the $i$th row of $D^{(1)}$ on that of $D^{(2)}$.
Throughout this paper, we let $D=(D^{(1)};D^{(2)})\in\R^{2n^2\times n^2}$,  $\rho(T)$ be the spectral radius of matrix $T$, and
$\mathcal{P}(\cdot)$ be the projection operator under Euclidean  norm.
The inner product of two vectors will be denoted by $\langle u,v\rangle$.
In the rest of this paper, we let $\|\cdot\| = \|\cdot\|_2$, and without misleading we abbreviate $\sum\nolimits_{i=1}^{n^2}$ as $\sum\nolimits_i$.
Additional notation will be introduced when it occurs.

The paper is organized as follows. In Section \ref{FTVCS}, we introduce our alternating minimization algorithm FTVCS and
study its convergence properties.  An accelerated scheme  of FTVCS is proposed in Section \ref{IADM} by incorporating an
inexact alternating direction technique. Numerical results in comparison with TwIST are presented in
Section \ref{NumericalResults}. Finally, we conclude the paper in
Section \ref{ConcludingRemarks}.

\section{Proposed algorithms}\label{Alg}

The task of this section is to construct our algorithm for solving \eqref{tv}. As is stated above, our interest in
this paper concentrates on CS encoding matrices, i.e., $A\in\R^{m\times n^2}$ with $m\ll n$.
The non-smoothness of TV causes the main difficulty. Similar as in \cite{YANG1}, we first
consider the approximation problem \eqref{atv} and then propose an inexact alternating direction method for
the solution of \eqref{tv}.

\subsection{Alternating minimization}\label{FTVCS}
The introduction of auxiliary variables $\vw$ in \eqref{atv} makes it easy to apply alternating minimization.
It is easy to see that, for fixed $u$,  the minimization of (\ref{atv}) with
respect to $\mathbf{w}$ reduces to the following two-dimensional problems
\begin{equation}\label{fixu}
\min_{\mathbf{w}_i\in\R^2}\|\mathbf{w}_i\|+\frac{\beta}{2}\|\mathbf{w}_i-D_iu\|^2,
\quad \ i=1,2,\ldots,n^2,
\end{equation}
for which the unique minimizers are given by the two-dimensional shrinkage formula
\begin{equation}\label{shrink}
\mathbf{w}_i=\max\left\{\|D_iu\|-\frac{1}{\beta}, \
0\right\}\frac{D_iu}{\|D_iu\|}, \ \ \ \ i=1,\dots,n^2,
\end{equation}
where the convention $0\cdot(0/0)=0$ is followed. On the other hand, for fixed $\mathbf{w}$,
the minimization of (\ref{atv}) with respect to $u$ is a least squares problem, and the corresponding normal equations
are given by
$$
\left(\sum\nolimits_{i}D_i^\top  D_i+\frac{\mu}{\beta}A^\top  A\right)u=\sum\nolimits_{i}D_i^\top  \mathbf{w}_i+\frac{\mu}{\beta}A^\top  f,
$$
or equivalently,
\begin{equation}\label{foru}
\left(D^\top  D+\frac{\mu}{\beta}A^\top  A\right)u=D^\top  w+\frac{\mu}{\beta}A^\top  f,
\end{equation}
where $w\in\R^{2n^2}$ is an reordering of $\vw_i$, $i=1,2,\ldots,n^2$.
It is well-known that, under the periodic boundary condition for
$u$, $D^\top  D$ is a block-circulant matrix and can be diagonalized
by two-dimensional fast Fourier transform (FFT). Unfortunately, the matrix $A^\top  A$ does not
have circulant structures for general CS encoding matrices. Therefore, the exact solution of
\eqref{foru} is expensive, which causes the main difficulty to apply alternating minimization directly.

To avoid solution of linear system of equations at each iteration, we
linearize $\frac{1}{2}\|Au-f\|^2$ at the current point $u^k$ and add a proximal term, resulting the following
approximation problem
\begin{equation}\label{line}
\min_{u,\mathbf{w}}\sum\nolimits_{i}\left(\|\mathbf{w}_i\|+\frac{\beta}{2}\|\mathbf{w}_i-D_iu\|^2\right)
+\mu\left(g_k^\top  (u-u^k)+\frac{1}{2\tau}\|u-u^k\|^2\right),
\end{equation}
where $g_k=A^\top  (Au^k-f)$ denotes the gradient of
$\frac{1}{2}\|Au-f\|^2$ at $u^k$, and $\tau>0$ is a parameter.
Clearly, problem (\ref{line}) is equivalent to
\begin{equation}\label{eqline}
\min_{u,\mathbf{w}}\sum\nolimits_{i}\left(\|\mathbf{w}_i\|+\frac{\beta}{2}\|\mathbf{w}_i-D_iu\|^2\right)
+\frac{\mu}{2\tau}\|u-(u^k-\tau g_k)\|^2.
\end{equation}
For fixed $w$ (or $\mathbf{w}$), the minimization of (\ref{eqline}) with respect to $u$
is equivalent to
\begin{equation}\label{eqnorm}
\left(D^\top  D+\frac{\mu}{\beta\tau}I\right)u=D^\top  w +\frac{\mu}{\beta\tau}(u^k-\tau
g_k),
\end{equation}
where we recall that $D=(D^{(1)}; D^{(2)})$. Under  the periodic boundary conditions for $u$,
the coefficient matrix in \eqref{eqnorm} can be diagonalized easily by FFT.
Consequently, the solution of \eqref{eqnorm} can be accomplished by two FFTs (including one inverse FFT).
To sum up, our alternating minimization algorithm, named fast total variation decoding from compressive sensing
measurements or FTVCS, is described  below.

\begin{algorithm}[FTVCS]\label{ftvdcsalg}
Input $f$, $A$ and $\mu,\beta,\tau>0$. Initialize $u^0=f$ and $k=0$.
\begin{tabbing}
(nr)ss\=ijkl\=bbb\=ccc\=ddd\= ee  \kill
\>{\bf While} ``not converged'', {\bf Do} \\
\>\>\>1) Compute $w^{k+1}$ according to (\ref{shrink}) for fixed $u=u^k$.\\
\>\>\>2) Compute $u^{k+1}$ according to (\ref{eqnorm}) for fixed $w = w^{k+1}$.\\
\>\>\>3) $k=k+1$.\\
\> {\bf End Do}
\end{tabbing}
\end{algorithm}

To establish the convergence of FTVCS, we need the following technical assumption.
\begin{assumption}\label{assum1}
$\mathcal{N}(A) \cap \mathcal{N}(D)=\{0\}$, where $\mathcal{N}(\cdot)$ represents the null space of a matrix.
\end{assumption}

Assumption \ref{assum1} is a quite loose condition and commonly used in the convergence analyses of
similar studies, see e.g., \cite{YANG1}. Under Assumption \ref{assum1}, we have the following convergence
results.

\begin{theorem}\label{conver}
Under Assumption \ref{assum1}, for any fixed $\beta>0$ and $0< \tau < 2 /\lambda_{\max}(A^\top  A)$, where
$\lambda_{\max}(A^\top A)$ denotes the spectral radius of $A^\top A$,
the sequence $\{(w^k,u^k)\}$ generated by Algorithm \ref{ftvdcsalg} from
any starting point $(w^0,u^0)$ converges to a solution
$(w^*,u^*)$ of \eqref{atv}.
\end{theorem}

\begin{theorem}\label{linear}
Suppose the sequence $\{(w^k,u^k)\}$ generated by Algorithm \ref{ftvdcsalg} converges to
$(w^*,u^*)$. Then, we have $\mathbf{w}_i^k=\mathbf{w}_i^*=0$,
$\forall \ i\in L$ after a finite number of iterations, where
$L=\{i:\|D_iu^*\|\leq 1/\beta\}$.
\end{theorem}

\begin{theorem}\label{qline}
Under the conditions of Theorem \ref{conver},
the sequence $\{u^k\}$ generated by Algorithm \ref{ftvdcsalg}
converges to $\{u^*\}$ $q$-linearly.
\end{theorem}

The proofs of Theorems \ref{conver}, \ref{linear} and \ref{qline} are given in Appendix \ref{AppendixA}.

\subsection{An accelerated scheme based on inexact alternating direction method}\label{IADM}
It is well-known that problem \eqref{atv} well approximates \eqref{tv} only when $\beta$ is sufficiently large.
However, it is generally difficult to determine theoretically how large a $\beta$ value must be
to attain a given accuracy. In this section, we present an inexact alternating direction method (ADM),
which converges to a solution of \eqref{tv} without requiring $\beta$ goes to
infinity.

First, we review briefly the idea of ADM pioneered in \cite{GLOWINSKI1,GLOWINSKI2}. The
classical ADM is designed to solve the following structure optimization problem:
\begin{equation}\label{adm}
\min_{y,z} \left\{\theta_1(y)+\theta_2(z):  Hy-z=0 \right\},
\end{equation}
where $\theta_1:\R^s\rightarrow\R$ and $\theta_2:\R^t\rightarrow\R$ are functions, and $H$ is a
$s\times t$ matrix. Given $z^k\in\R^t$ and $p^k\in\R^s$, the ADM iterates as follows
\begin{eqnarray}\label{ADM-G}
y^{k+1}&\leftarrow&\mbox{arg}\min_y  \theta_1(y)- (p^k)^\top (Hy-z^k)+\frac{\sigma}{2}\|Hy-z^k\|^2,\\
z^{k+1}&\leftarrow&\mbox{arg}\min_z \theta_2(z)- (p^k)^\top(Hy^{k+1}-z)+\frac{\sigma}{2}\|Hy^{k+1}-z\|^2,\\
p^{k+1} &\leftarrow& p^k-\sigma(Hy^{k+1}-z^{k+1}),
\end{eqnarray}
where $\sigma>0$ is a parameter. In \eqref{ADM-G},  $p^k$ is the Lagrangian multiplier and $\sigma$ severs as
a penalty parameter. It can be shown that, under quite reasonable assumption, \eqref{ADM-G} converges to
a solution of \eqref{adm} for any fixed $\sigma>0$, see  \cite{GLOWINSKI1,GLOWINSKI2}.

We now consider the model (\ref{tv}) in its equivalent form
\begin{equation}\label{eqftvdcs}
\min_{u,w}\left\{ \sum\nolimits_{i}\|\mathbf{w}_i\|+\frac{\mu}{2}\|Au-f\|^2: \mathbf{w}_i=D_i u, \forall i\right\}.
\end{equation}
The augmented Lagrangian problem of \eqref{eqftvdcs} is given by
\begin{equation}\label{ALP}
\min_{u,w}\sum\nolimits_{i}\left(\|\mathbf{w}_i\| - \lambda_i^\top(\mathbf{w}_i-D_i u) +
\frac{\beta}{2}\|\mathbf{w}_i-D_i u\|^2\right)+\frac{\mu}{2}\|Au-f\|^2,
\end{equation}
where, for each $i$, $\lambda_i\in\R^2$ is the Lagrangian multiplier attached to $\mathbf{w}_i=D_i u$.
Inspired by the ADM iterations, for given $(u^k,w^k,\lambda^k)$,
we obtain the next triplet $(u^{k+1},w^{k+1},\lambda^{k+1})$ as follows.
First, for fixed $u^k$ and $\lambda^k$, the minimization of \eqref{ALP} with respect to $\vw$ is equivalent to
\begin{eqnarray*}
\min_{\mathbf{w}_i\in\R^2} \|\mathbf{w}_i\|+\frac{\beta}{2}\|\mathbf{w}_i-(D_iu^k- \lambda_i^k/\beta)\|^2, \; i=1,2,\ldots, n^2,
\end{eqnarray*}
the solutions of which are given by
\begin{equation}\label{iadmw}
\mathbf{w}_i^{k+1}=\max\left\{\|D_iu^k-\lambda_i^k/\beta\| - \frac{1}{\beta}, 0\right\}\frac{D_iu^k-
\lambda_i^k/\beta}{\|D_iu^k-\lambda_i^k/\beta\|},
 \ \ \ i=1,2,\ldots,n^2.
\end{equation}
Second, for fixed $w^{k+1}$, $u^k$ and $\lambda^k$,
the minimization of \eqref{ALP} with respect to $u$ is
approximated by linearizing $\frac{1}{2}\|Au-f\|^2$ and adding a
proximal term as in \eqref{line}, resulting the following problem
\begin{equation}\label{linetv}
\min_{u} \sum\nolimits_{i} \left(- (\lambda_i^k)^\top (\mathbf{w}_i^{k+1}-D_iu) + \frac{\beta}{2}\|\mathbf{w}_i^{k+1}-D_iu\|^2\right)
+\frac{\mu}{2\tau}\|u-(u^k-\tau
g_k)\|^2,
\end{equation}
where $g_k$ is defined in (\ref{line}). It is easy to show that the normal equations of \eqref{linetv} are of the form
\begin{equation}\label{eqnorm11}
\left(D^\top  D+\frac{\mu}{\beta\tau}I\right)u=D^\top  (w^{k+1}-  \lambda^k/\beta)+\frac{\mu}{\beta\tau}(u^k-\tau
g_k).
\end{equation}
Under the periodic boundary conditions, the exact solution of \eqref{eqnorm11} can be attained
by two FFTs. Finally, $\lambda$ is updated via
\begin{equation}\label{udlamb}
\lambda^{k+1}=\lambda^k- \beta (w^{k+1}-Du^{k+1}),
\end{equation}
We note that the linearization technique makes the
$u$-subproblem of \eqref{ALP} is solved inexactly. Therefore, we name  the above iterative framework
as an inexact ADM or IADM, which is summarized below.

\begin{algorithm}[IADM]\label{iadmalg}
Input $f$, $A$ and $\mu,\beta,\tau>0$. Initialize $u^0=f$ and $k=0$.
\begin{tabbing}
(nr)ss\=ijkl\=bbb\=ccc\=ddd\= ee  \kill
\>{\bf While} ``not converged'', {\bf Do} \\
\>\>\>1) Compute $w^{k+1}$ according to (\ref{iadmw}) for fixed $\lambda=\lambda^k$ and $u=u^k$.\\
\>\>\>2) Compute $u^{k+1}$ according to (\ref{eqnorm11}) for fixed $\lambda=\lambda^k$ and $w = w^{k+1}$.\\
\>\>\>3) Update $\lambda$ via \eqref{udlamb} and set $k=k+1$.\\
\> {\bf End Do}
\end{tabbing}
\end{algorithm}

We have the following convergence results for Algorithm \ref{iadmalg}.
\begin{theorem}\label{augconv}
Under Assumption \ref{assum1}, the sequence $\{(w^k,u^k)\}$ generated by Algorithm \ref{iadmalg} from any starting point
$(w^0,u^0)$ converges to a solution of \eqref{eqftvdcs}.
\end{theorem}

A closer examination shows that Algorithm \ref{iadmalg} is related to
the proximal ADM of He et al. \cite{BSHE} for solving monotone variational inequalities.
Hence, the global convergence is followed directly, see Appendix \ref{AppendixB} for details.

\section{Numerical experiments} \label{NumericalResults}

In this section, we present numerical results to illustrate
the feasibility and efficiency of FTVCS and its accelerated variant IAMD. All experiments were
accomplished in Matlab 2009a running on a PC ({\ttfamily Intel Pentium(R) 4, 1.6 GHz,
1.0GB SDRAM}) with Windows XP operating system. As usual, we measure
the quality of reconstruction by relative error  to the
original image $\bar{u}$, i.e.,
$$
RE= \frac{\|u-\bar{u}\|}{\|\bar{u}\|} \times 100\%.
$$
In the following, we first present primary  experimental results
to show the feasibility of both algorithms,
and then compare both algorithms with TwIST --- a state-of-the-art algorithm
for solving \eqref{tv}, to demonstrate their efficiency.

\subsection{Test on FTVCS and IADM}

\begin{figure}[htbp]
\vspace{0.2cm} \centering
\includegraphics[scale = .55]{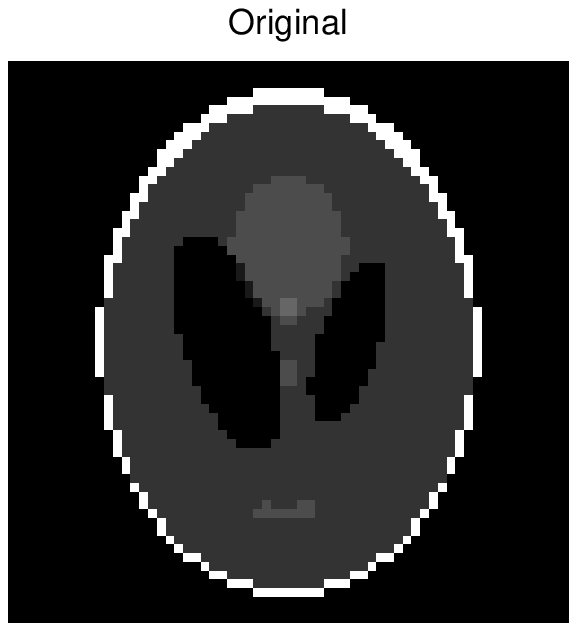}\hspace{.2cm}
\includegraphics[scale = .55]{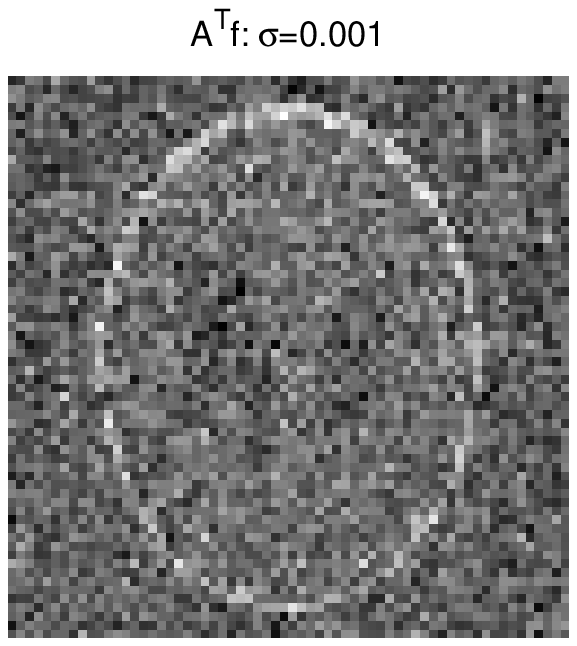}\hspace{.2cm}
\includegraphics[scale = .55]{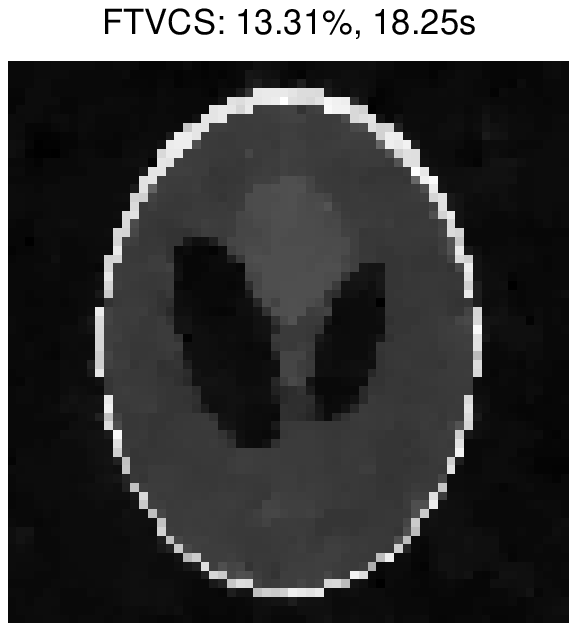}\hspace{.2cm}
\includegraphics[scale = .55]{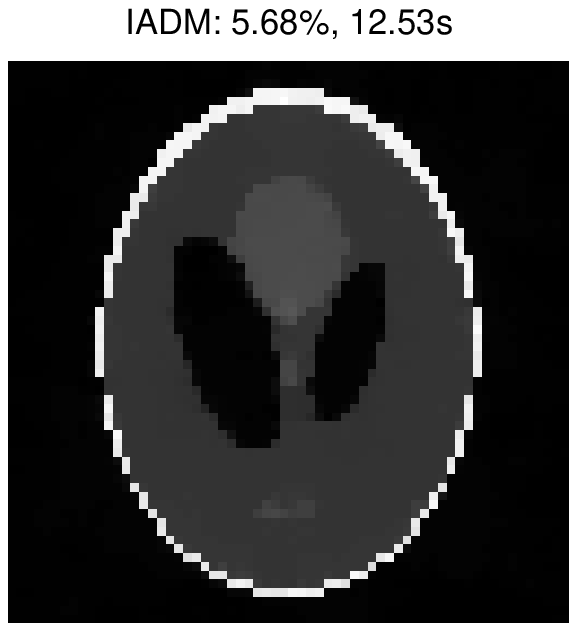}\hspace{.2cm}
\vspace{0cm} \caption{{\footnotesize Reconstruction results of FTVCS and
IAMD. Original (left, $64\times 64$); Initial guess (middle
left); Recovered by FTVCS (middle right, RE=$13.31\%$,
CPU time $18.25$s); Recovered by IAMD (right,
RE=$5.68\%$, CPU time $12.53$s)}}. \label{figure1} \vspace{-0.5cm}
\end{figure}

In the first experiment, we present reconstruction results of both algorithms
to illustrate their feasibility for solving \eqref{tv}. We used a random
matrix with independent identical distributed Gaussian entries as CS encoding matrix
and tested the Shepp-Logan phantom image, which has been widely used in simulations for
TV models.  Due to storage limitations, we tested the image size $64\times 64$.
The sample ratio in this test is $30\%$, which are selected uniformly at random.
Besides, we added Gaussian noise of zero mean and standard deviation
$\sigma=0.001$.  Similar as in FTVd \cite{YANG1}, we implemented FTVCS
with a continuation scheme on $\beta$ to speed up convergence.
Specifically, we tested the $\beta$-sequence $\{2^4, 2^5,2^6,2^7\}$ and used
the warm-start technique. In IADM, the value of $\beta$ is fixed to be $8$.
In both algorithms, the weighting parameter $\mu$ was set to be $200$.
Both algorithms were terminated when the relative change between successive iterates fell
below $10^{-3}$, i.e.,
\begin{equation}\label{stop}
\|u_k-u_{k-1}\|\leq 10^{-3}\|u_{k-1}\|.
\end{equation}
The original image, the initial guess, and the reconstructed ones by
both algorithms are listed in Figure \ref{figure1}.

\begin{figure}[htbp]
\vspace{-0cm}\centering
\includegraphics[scale = 0.45]{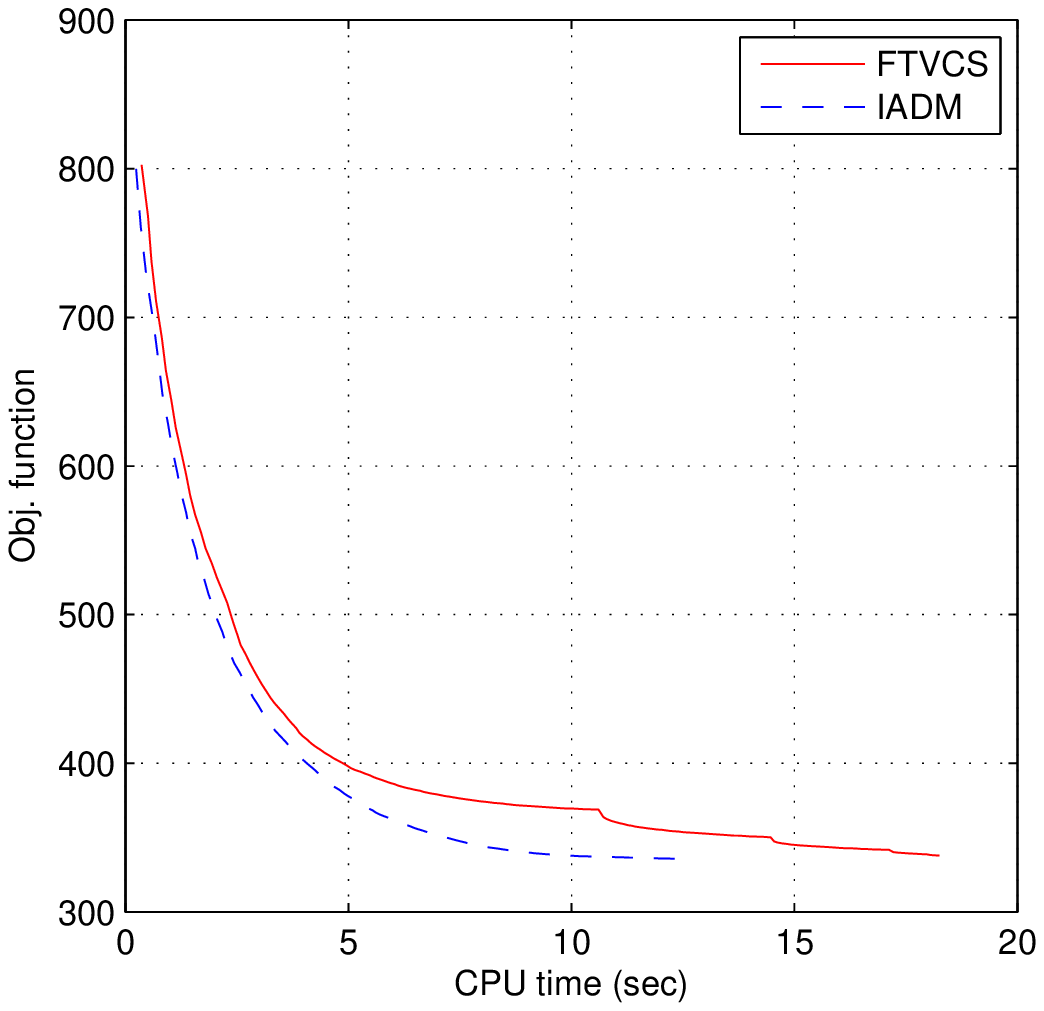}\hspace{-0.0cm}
\includegraphics[scale = 0.45]{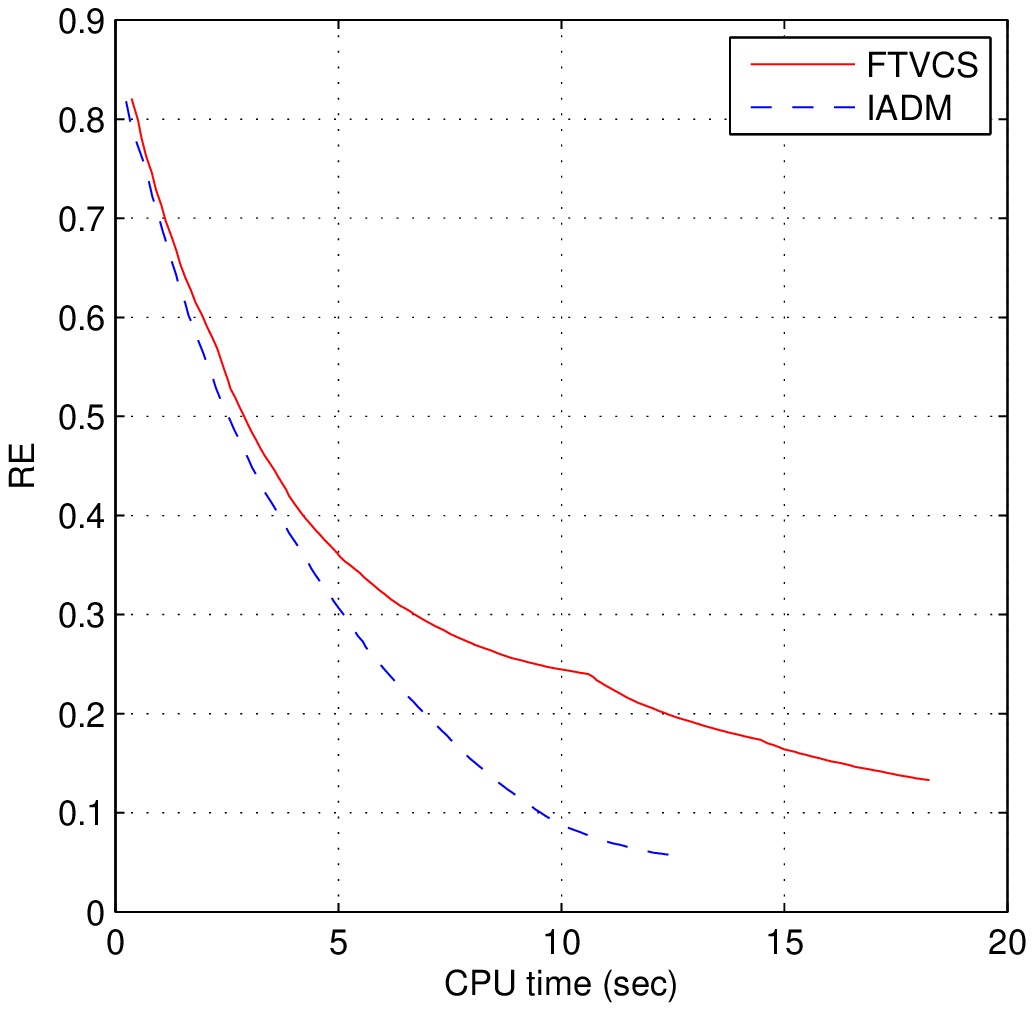}
\vspace{-0.5cm} \caption{{\footnotesize Convergence
behavior of FTVCS and IADM. Left: objective function; Right: relative error.
In both plots, the horizontal axes denote CPU time in seconds.}}\label{figure2}
\end{figure}

As is shown in Figure \ref{figure1}, both algorithms perform favorably and produce faithful recovery results
in a few seconds. We note that the per-iteration cost of both algorithms
is one shrinkage operation, two matrix-vector multiplications and two FFTs.
The results also indicate that the inexact ADM approach described in Algorithm \ref{iadmalg}
is indeed more efficient than the penalty approach FTVCS described in
Algorithm \ref{ftvdcsalg} in the sense that better recovery results were obtained in less CUP seconds.
To closely examine the convergence behavior of both algorithms, we present in Figure \ref{figure2}
the decreasing of objective function values and  relative errors as CPU time proceeded.
It is clear from Figure \ref{figure2} that both algorithms generated decreasing sequences of function values.
From the right-hand plot, IADM achieved a solution of lower relative error. In both plots,
the curves of IADM fall bellow those of FTVCS throughout the whole iteration process.

\subsection{Comparison with TwIST}
In this subsection, we present extensive numerical results to compare IADM with TwIST \cite{TWIST} ---
a two-step iterative shringkage/thresholding algorithm for solving a class of
optimization problems arising from image restoration, reconstruction and linear inverse
problems\footnote{The Matlab code of TwIST can be obtained from
\url{http://www.lx.it.pt/~bioucas/TwIST/TwIST.htm}
}. Specifically, TwIST is designed to
solve
\begin{equation}\label{lip}
\min_u \mathcal{J}(u)+\frac\mu2\|Au-b\|^2,
\end{equation}
where $\mathcal{J}(\cdot)$ is a general regularizer, which can be
either the $\ell_1$-norm or the TV semi-norm, as well as others.
In the comparison, we used  partial discrete cosine transform (DCT)
matrix as CS encoder, i.e., the $m$ rows of $A$ were chosen uniformly
at random from the $n\times n$ DCT matrix. Since the DCT
matrix is implicity stored as fast transforms, this enables us
to test larger images.
We used the default parametric settings for TwIST and terminated it
as the relative change in objective function values fell below $tol=10^{-3}$.
The parameters in IADM were set as follows: $\tau=1.9$ and $\beta=2^6$.
To obtain higher quality images, we used more stringent stopping tolerance and
terminated IADM when $\|u_k-u_{k-1}\|\leq 5\times10^{-5}\|u_{k-1}\|$ was satisfied.

We first compared IADM with TwIST using the Shepp-Logan phantom benchmark image of
size $128\times 128$. We randomly selected $30\%$ DCT coefficients and added
Gaussian noise of mean zero and standard deviation $0.001$.
Table \ref{table1} reports the detailed results of both algorithms
for different values of $\mu$, where RE, Obj, Iter and Time represent, respectively,
the relative error of the reconstructed image to the original one, the final
objective function value, the number of iterations, and the consumed CPU time in
seconds.

It can be seen from Table \ref{table1} that, for both algorithms,
the number of iterations becomes larger and larger as $\mu$ increases, and as a result
longer CPU time is consumed. For larger $\mu$, the performance of both algorithms deteriorates, while
the resulting relative errors were not improved. For $\mu$ between $500$ and $7000$,
IADM always obtained comparable or higher recovery quality than TwIST.
For $\mu$ between $500$ and $1000$, IADM is also faster than TwIST.
In terms of final function values, IADM obtained slightly smaller ones than those of TwIST.

\begin{table}
\begin{center}
\scriptsize{  \caption{Comparison results of IADM and TwIST
with different $\mu$.}
\begin{tabular}{l|rrrr|rrrr}
\hline\hline \multicolumn{1}{c|}{ } & \multicolumn{4}{c|}{TwIST } &
\multicolumn{4}{c}{ IADM} \\ \hline
$\mu$ & RE  &  Obj  &  Iter  &  Time  & RE  &  Obj  &  Iter  &  Time\\
\hline
  100 & 3.79\%  &      736.10  &  40  &  26.51  & 4.80\%  &  721.95  &  140  &   12.73\\
  500 & 3.74\%  &      804.20  &  55  &  29.23  & 3.37\%  &  786.99  &  219  &   19.26\\
  600 & 3.93\%  &      810.80  &  56  &  29.02  & 3.43\%  &  794.20  &  247  &   21.47\\
  700 & 3.97\%  &      814.32  &  62  &  33.27  & 3.56\%  &  798.19  &  271  &   23.34\\
  800 & 3.90\%  &      815.41  &  69  &  33.17  & 3.59\%  &  800.68  &  297  &   23.16\\
  900 & 4.06\%  &      818.17  &  65  &  31.92  & 3.51\%  &  803.09  &  310  &   26.97\\
 1000 & 4.71\%  &      817.52  &  71  &  29.72  & 3.81\%  &  803.00  &  354  &   25.86\\
 2000 & 4.45\%  &      831.86  &  91  &  43.41  & 3.99\%  &  817.08  &  602  &   53.22\\
 3000 & 5.40\%  &      926.04  &  86  &  34.55  & 4.54\%  &  816.49  &  845  &   75.22\\
 4000 & 4.40\%  &      830.97  & 122  &  47.23  & 4.23\%  &  817.81  & 1066  &   84.11\\
 5000 & 4.54\%  &      873.51  & 146  &  61.06  & 4.23\%  &  821.52  & 1303  &  116.03\\
 6000 & 4.51\%  &      858.31  & 182  &  77.19  & 4.33\%  &  822.30  & 1581  &  132.27\\
 7000 & 4.61\%  &      852.42  & 200  &  89.05  & 4.39\%  &  821.55  & 1742  &  157.08\\
 8000 & 4.10\%  &      832.80  & 288  & 139.17  & 4.42\%  &  821.31  & 1929  &  177.05\\
 9000 & 4.22\%  &      830.45  & 315  & 135.16  & 4.59\%  &  819.54  & 2165  &  170.30\\
10000 & 4.10\%  &      828.36  & 336  & 156.98  & 4.47\%  &  817.46  & 2405  &  245.77\\
\hline\hline
\end{tabular}
}\label{table1}
\end{center}
\end{table}

\begin{figure}[htbp]
\vspace{-0.0cm}  \centering
\includegraphics[scale = .37]{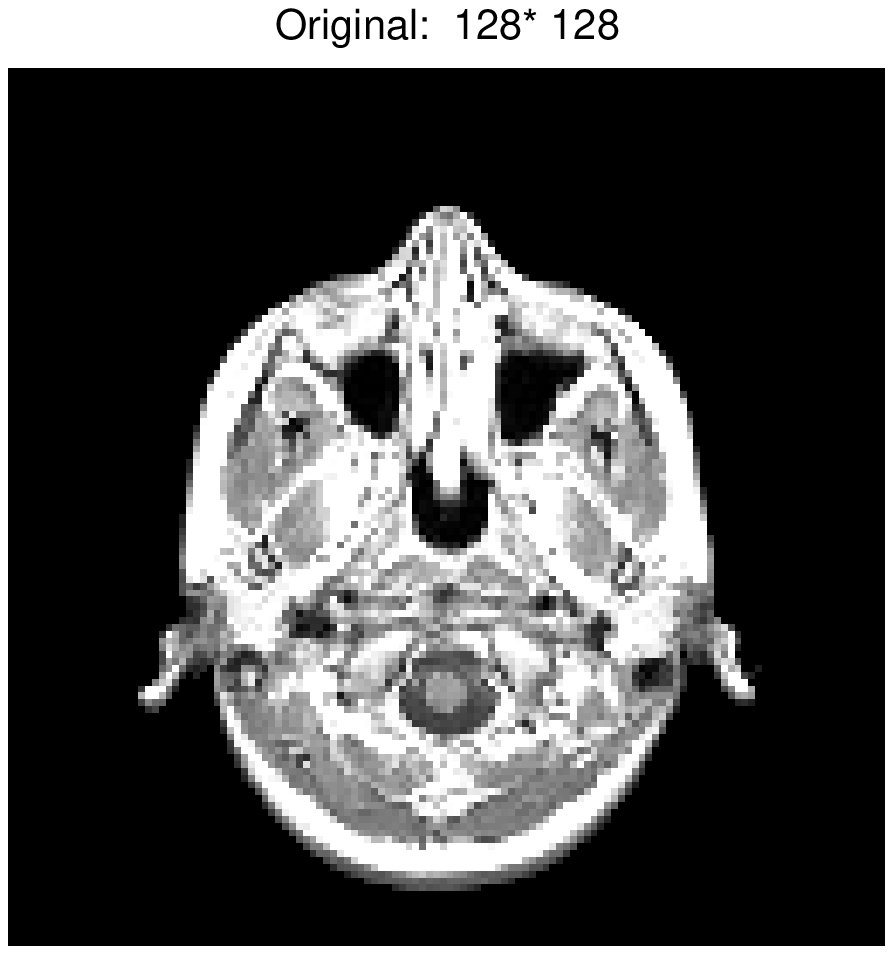}\hspace{.2cm}
\includegraphics[scale = .37]{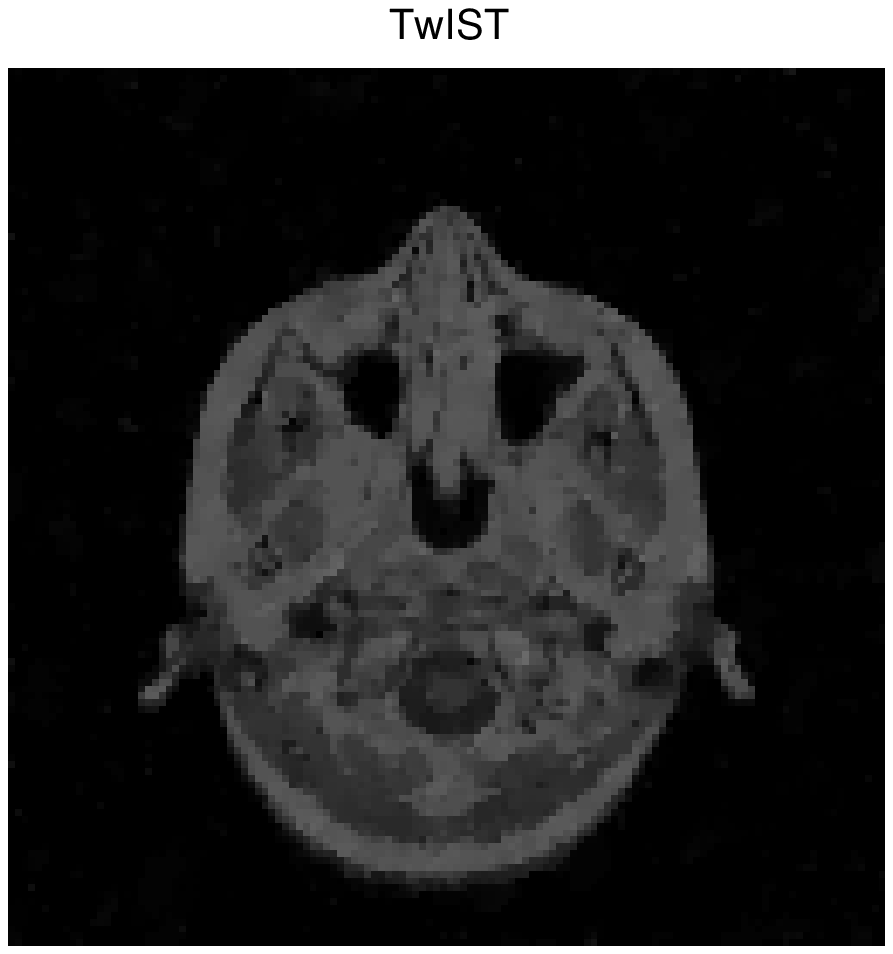}\hspace{.2cm}
\includegraphics[scale = .37]{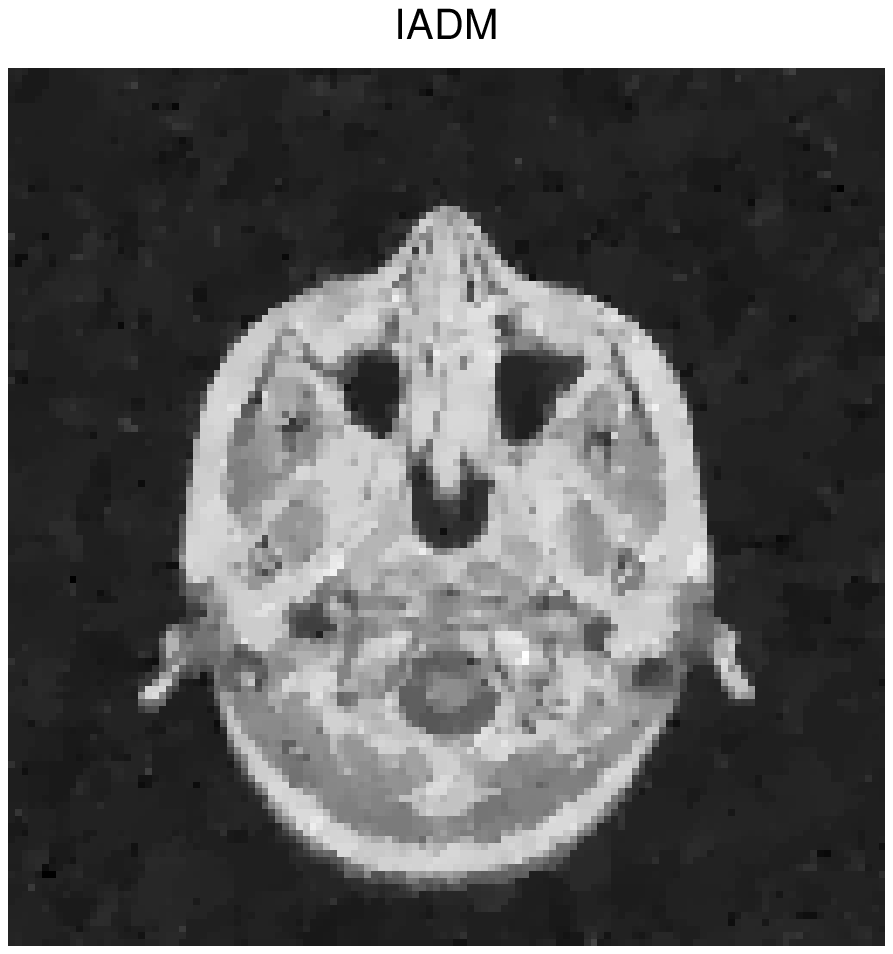}\hspace{.2cm}\\
\vspace{0.5cm}
\includegraphics[scale = .29]{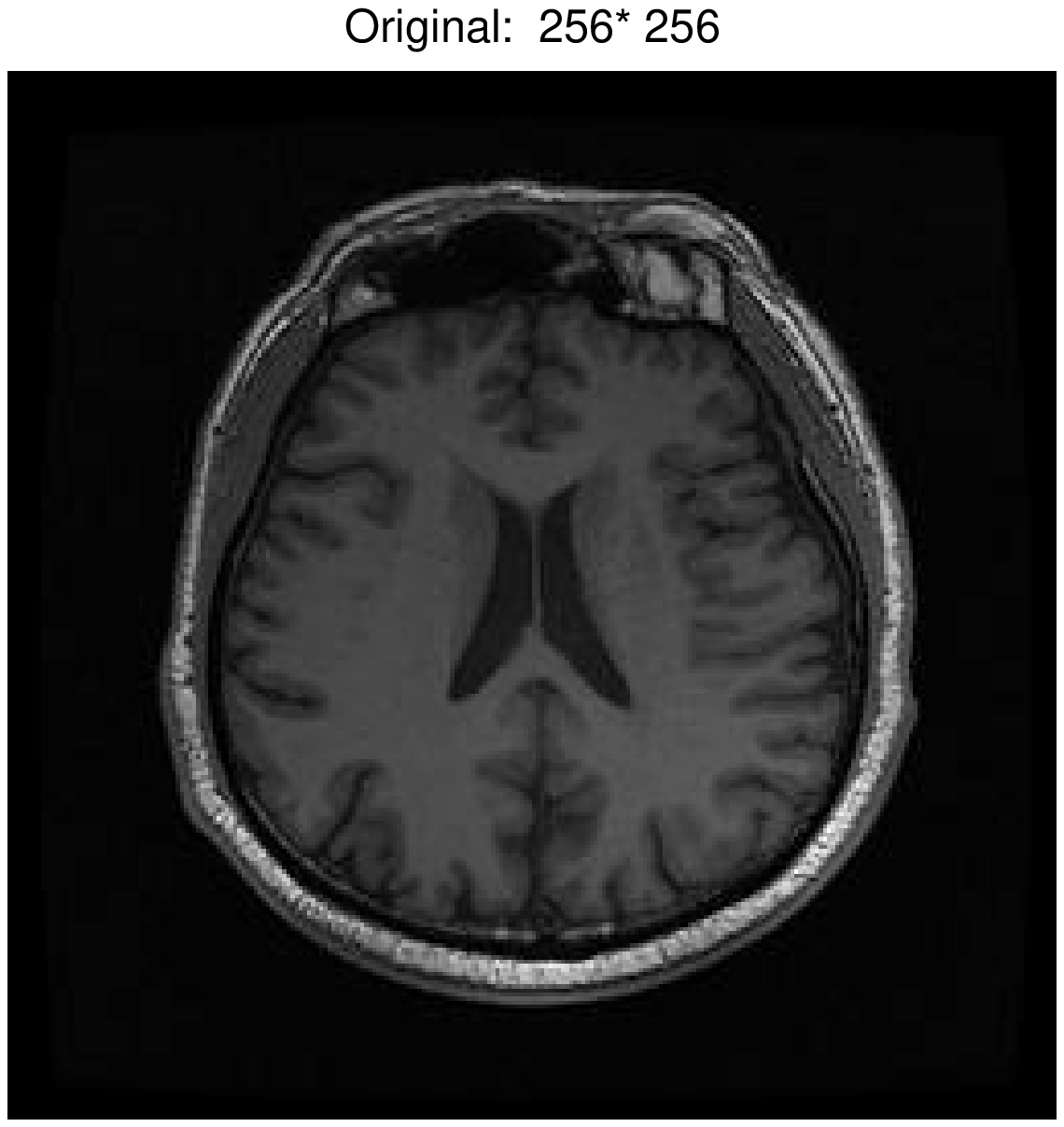}\hspace{.2cm}
\includegraphics[scale = .29]{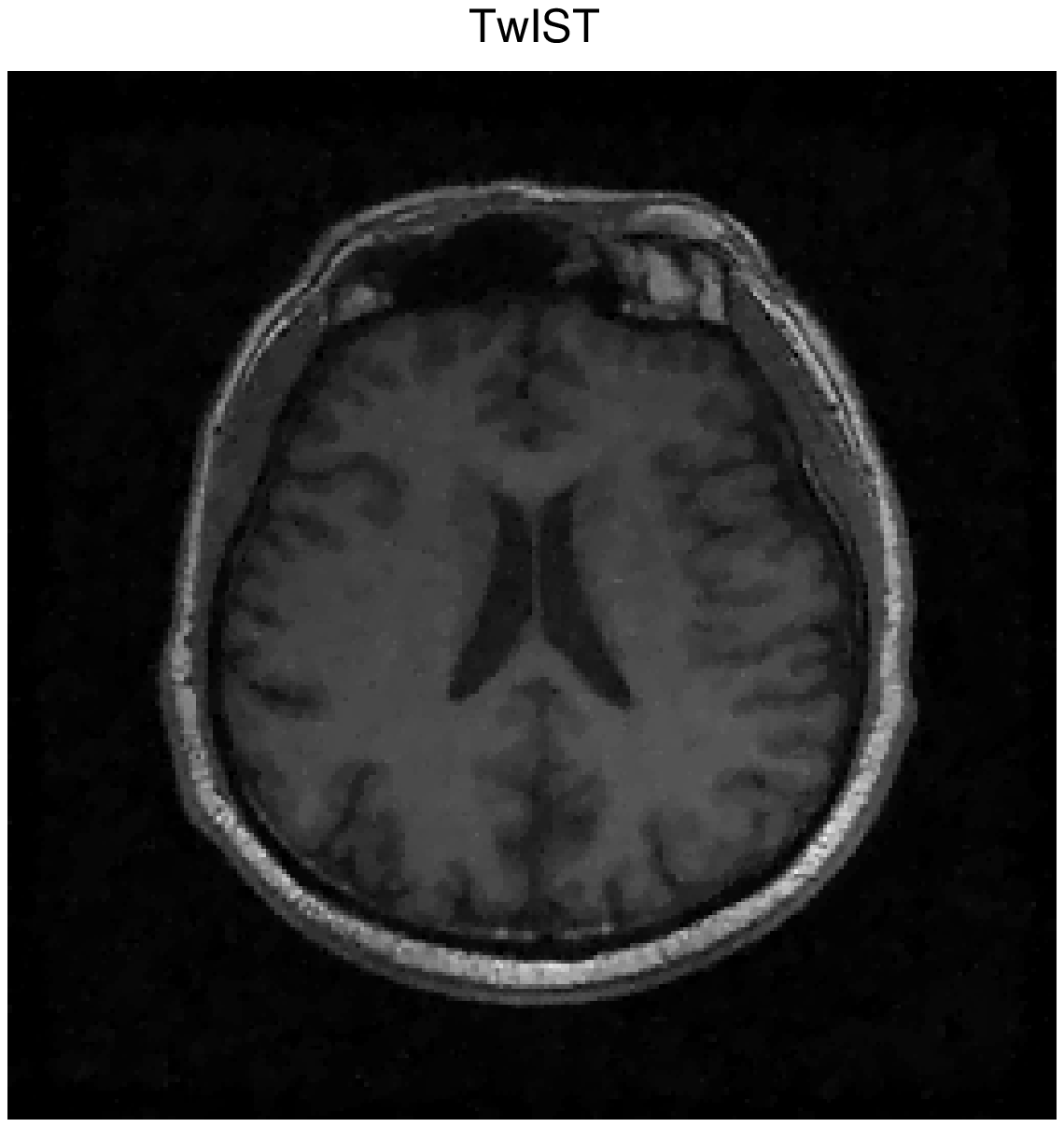}\hspace{.2cm}
\includegraphics[scale = .29]{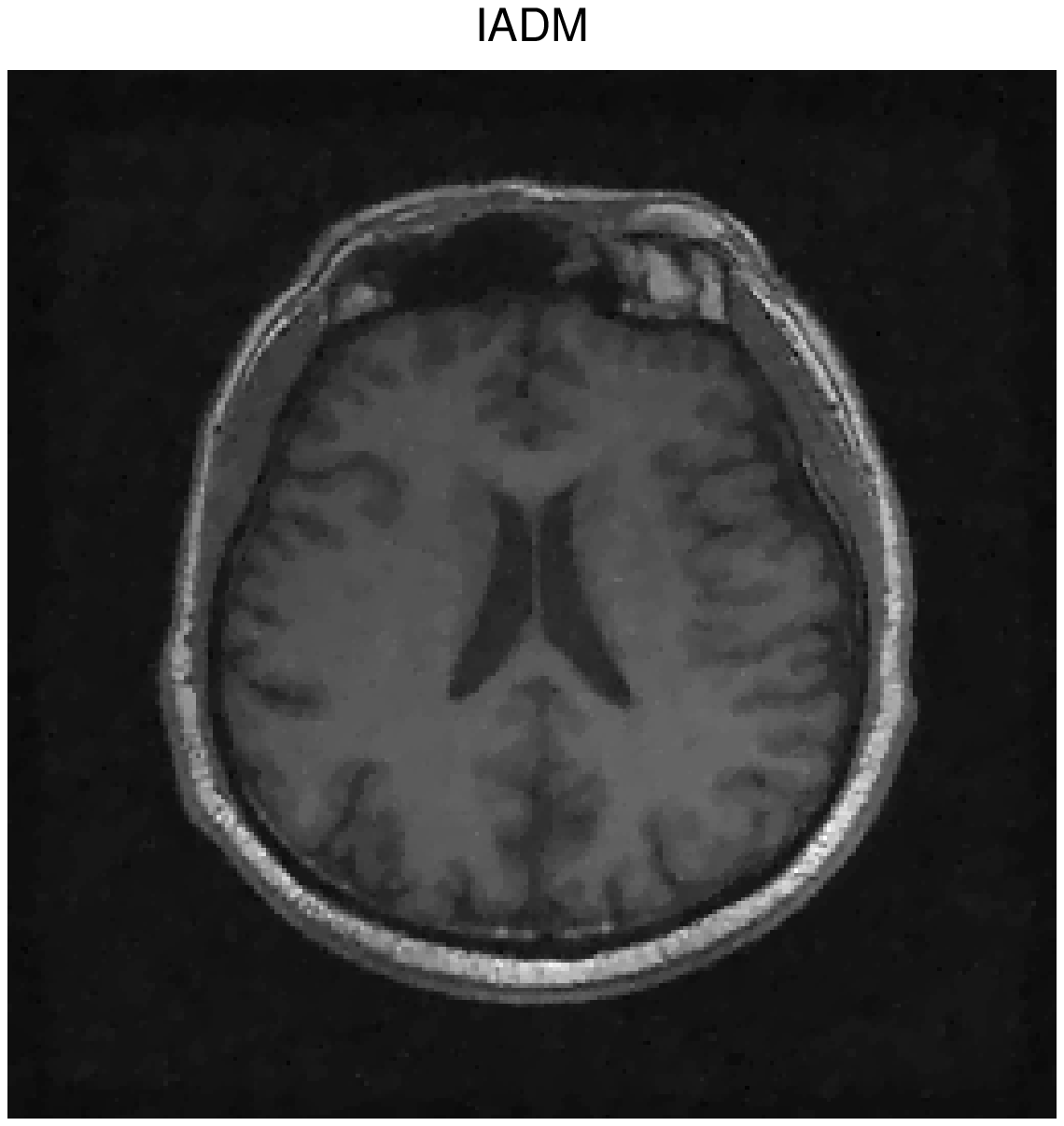}\hspace{.2cm}\\
\vspace{0.5cm}
\includegraphics[scale = .29]{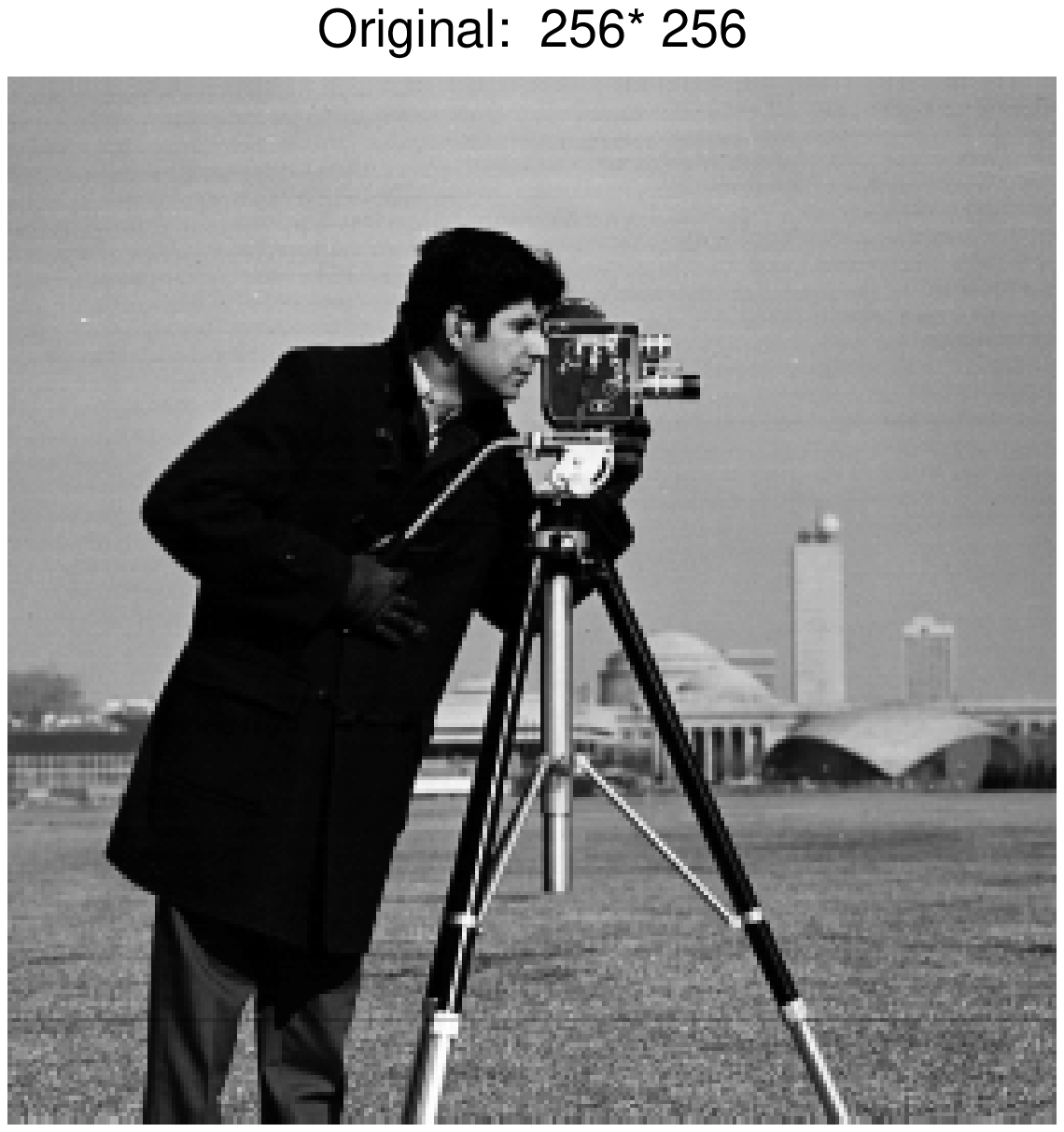}\hspace{.2cm}
\includegraphics[scale = .29]{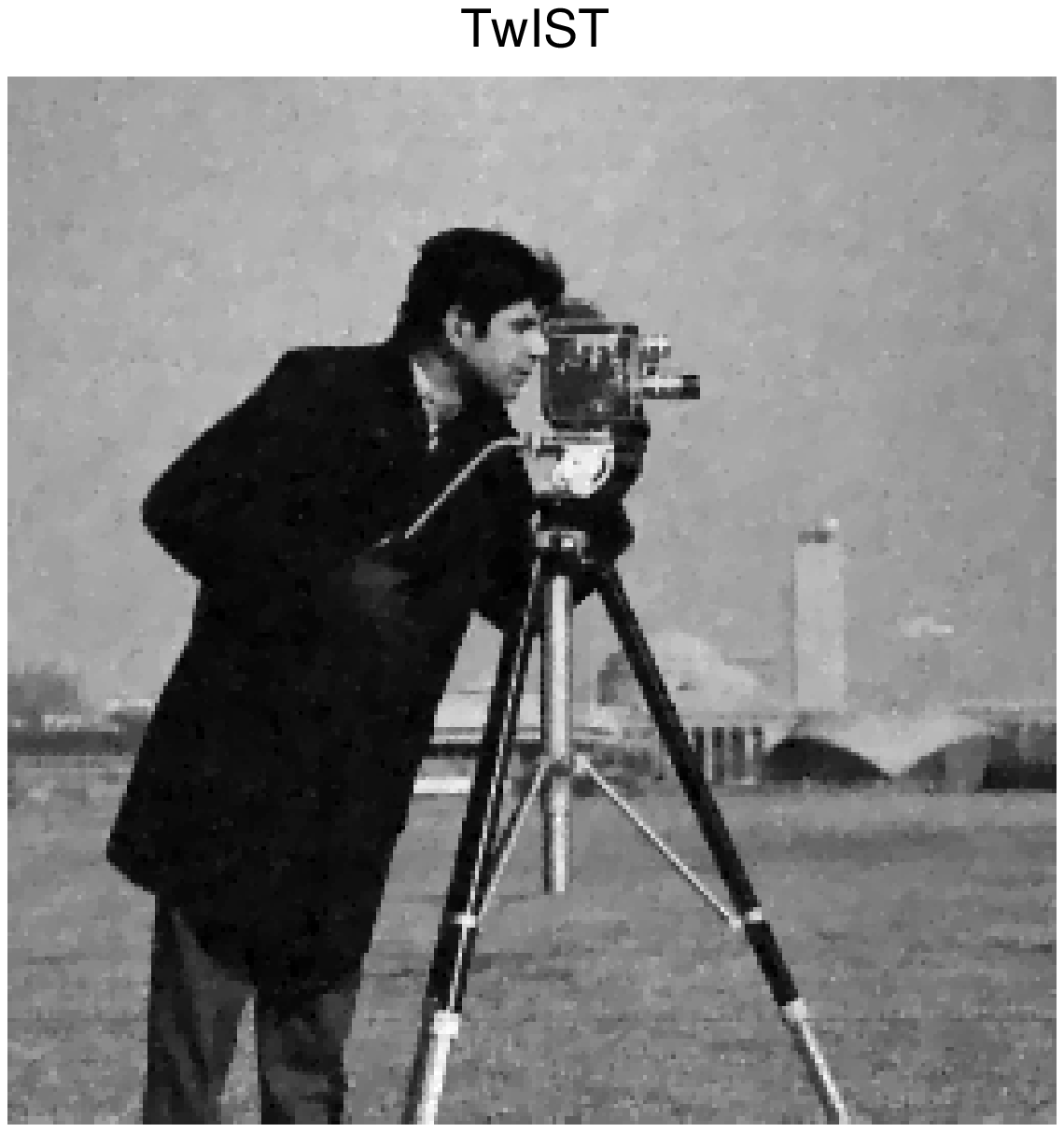}\hspace{.2cm}
\includegraphics[scale = .29]{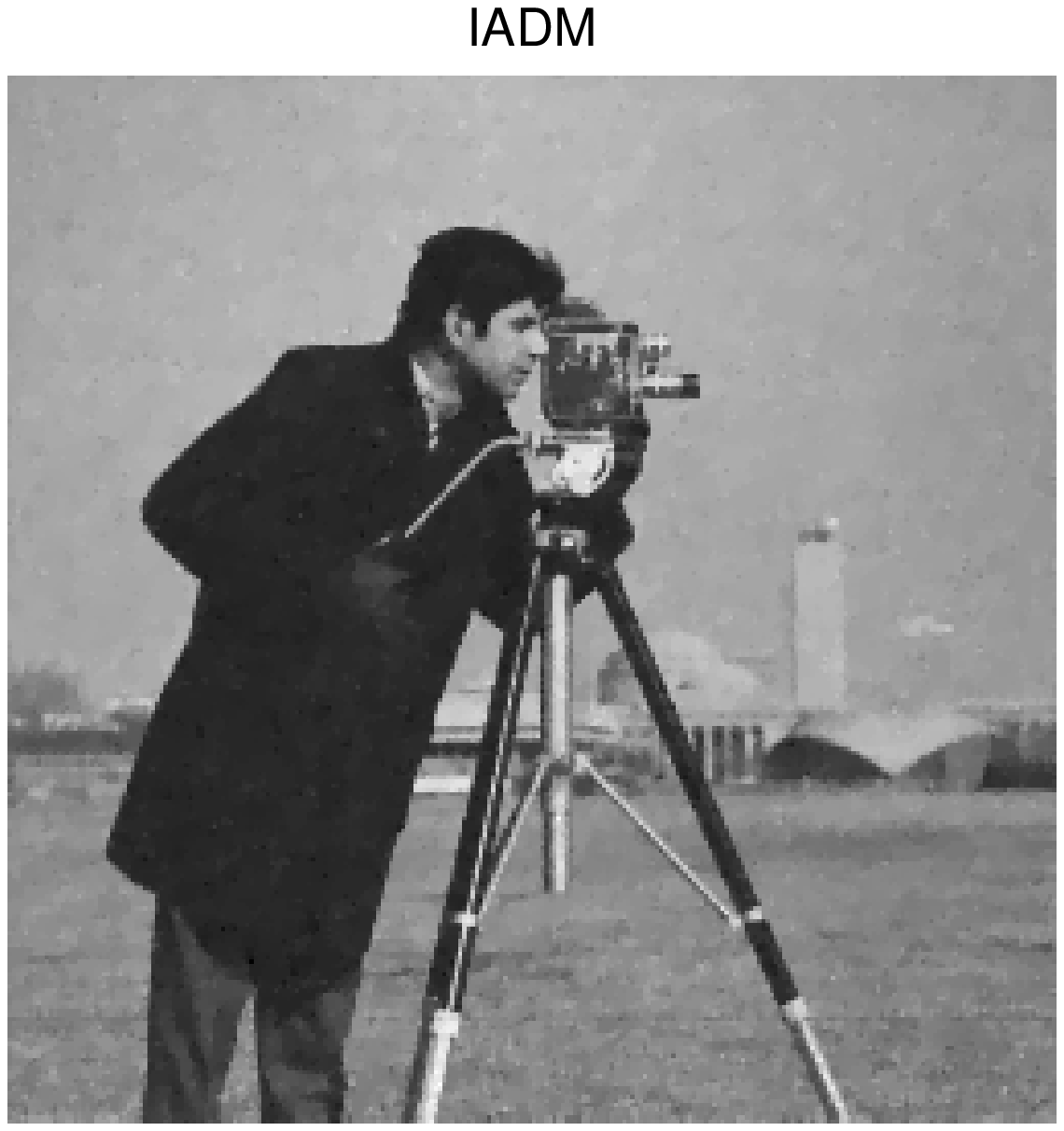}\hspace{.2cm}\\
\vspace{0.5cm}
\includegraphics[scale = .29]{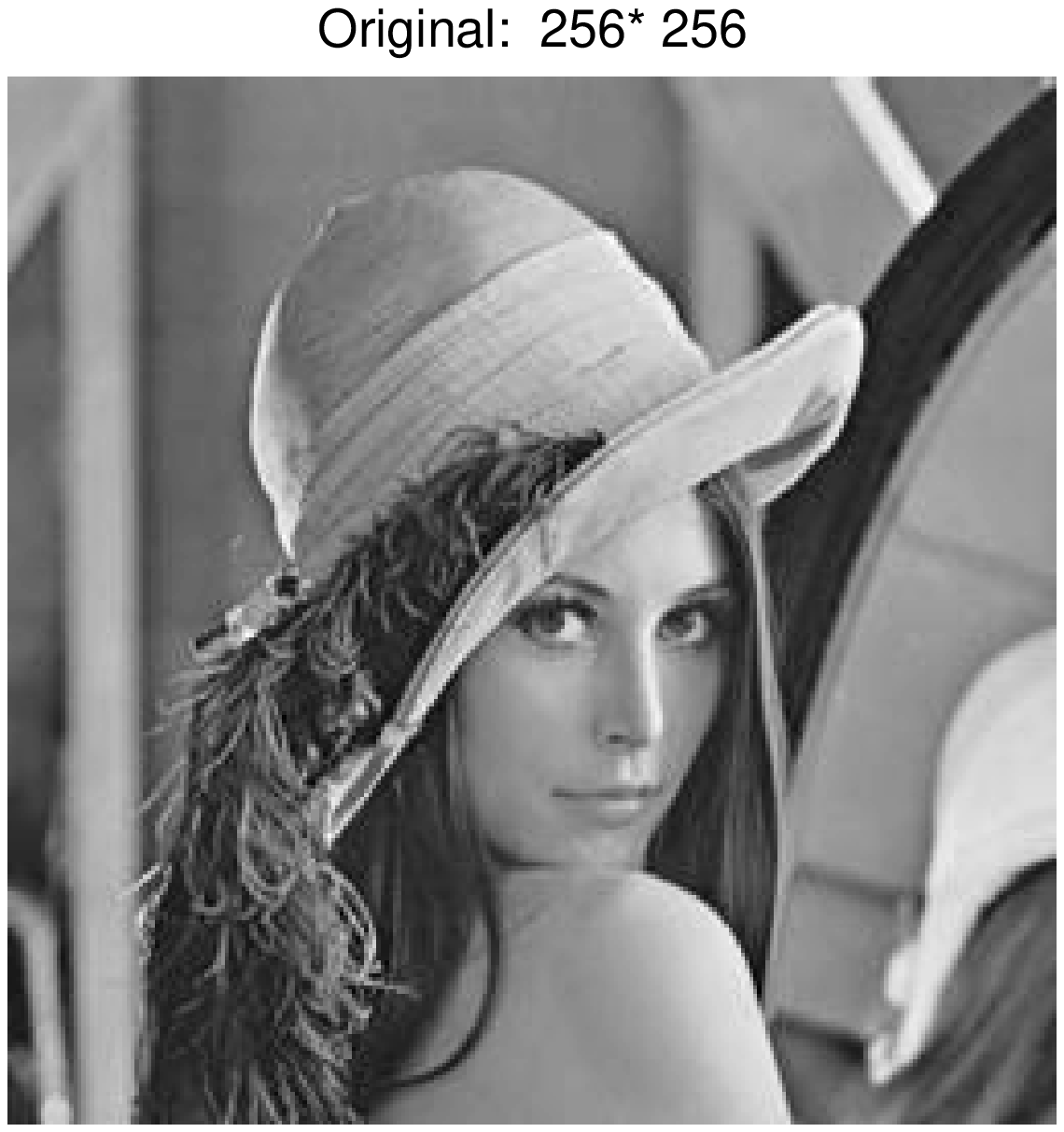}\hspace{.2cm}
\includegraphics[scale = .29]{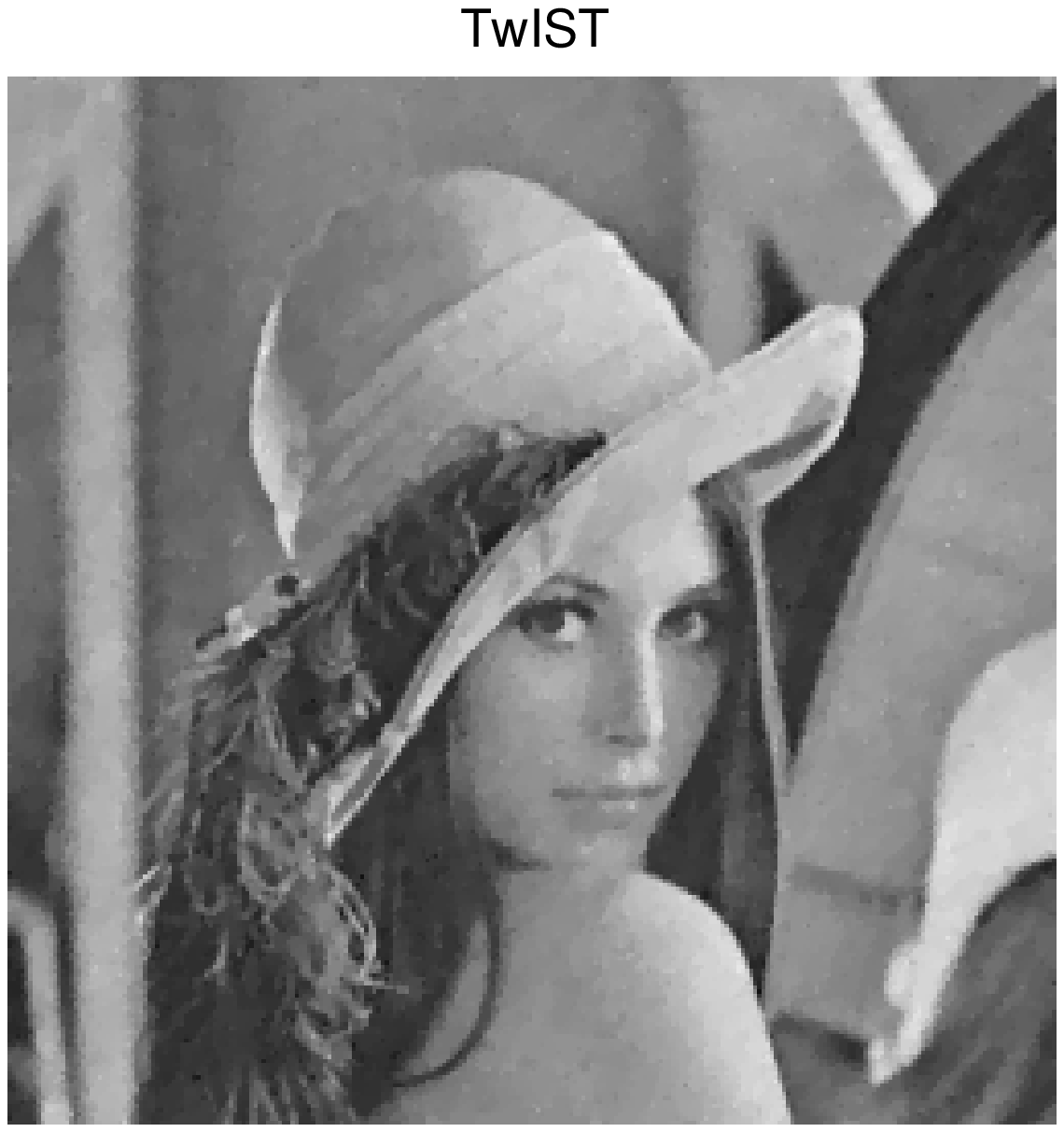}\hspace{.2cm}
\includegraphics[scale = .29]{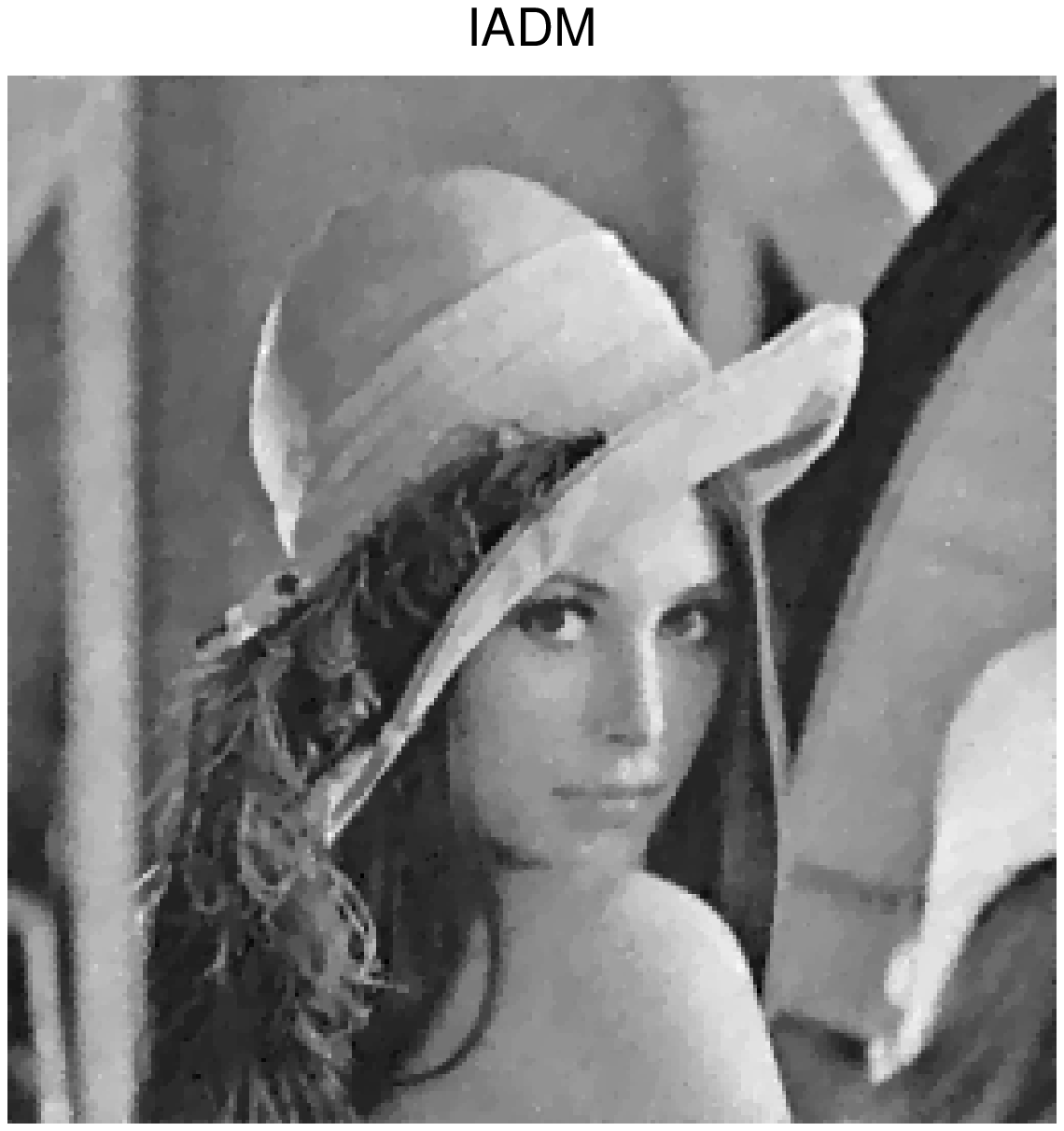}\hspace{.2cm}\\
\vspace{0.5cm}
\includegraphics[scale = .272]{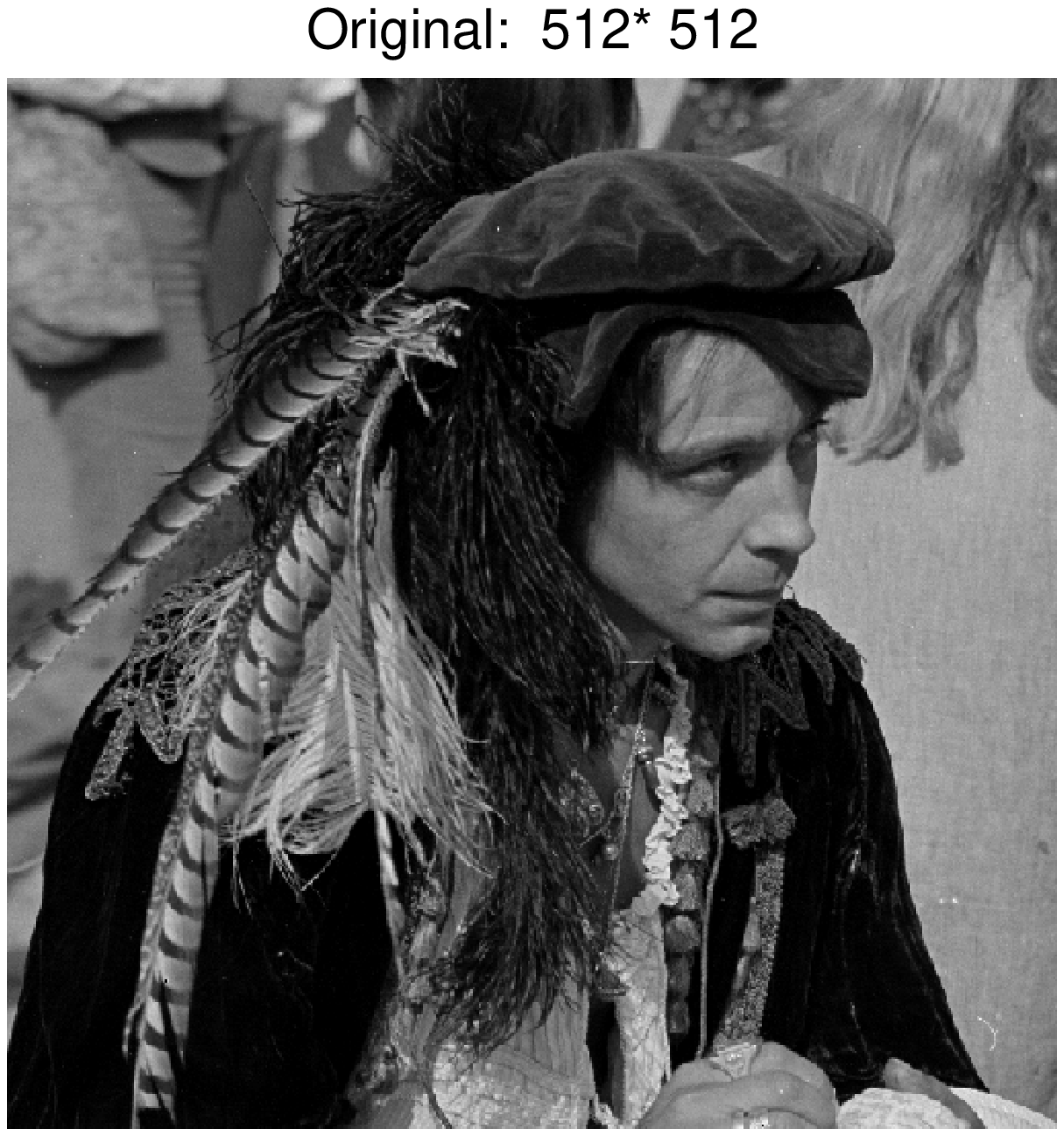}\hspace{.2cm}
\includegraphics[scale = .272]{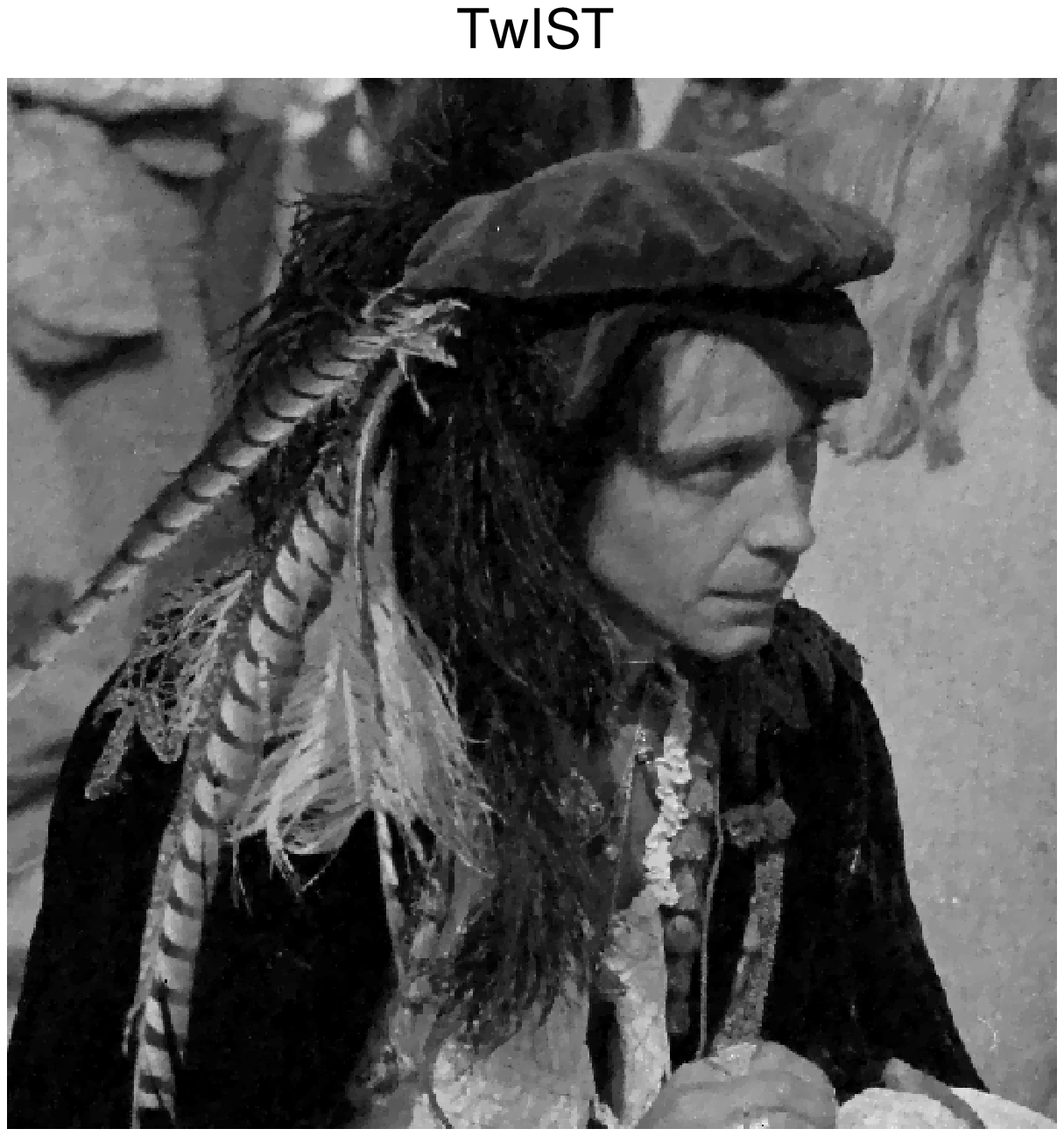}\hspace{.2cm}
\includegraphics[scale = .272]{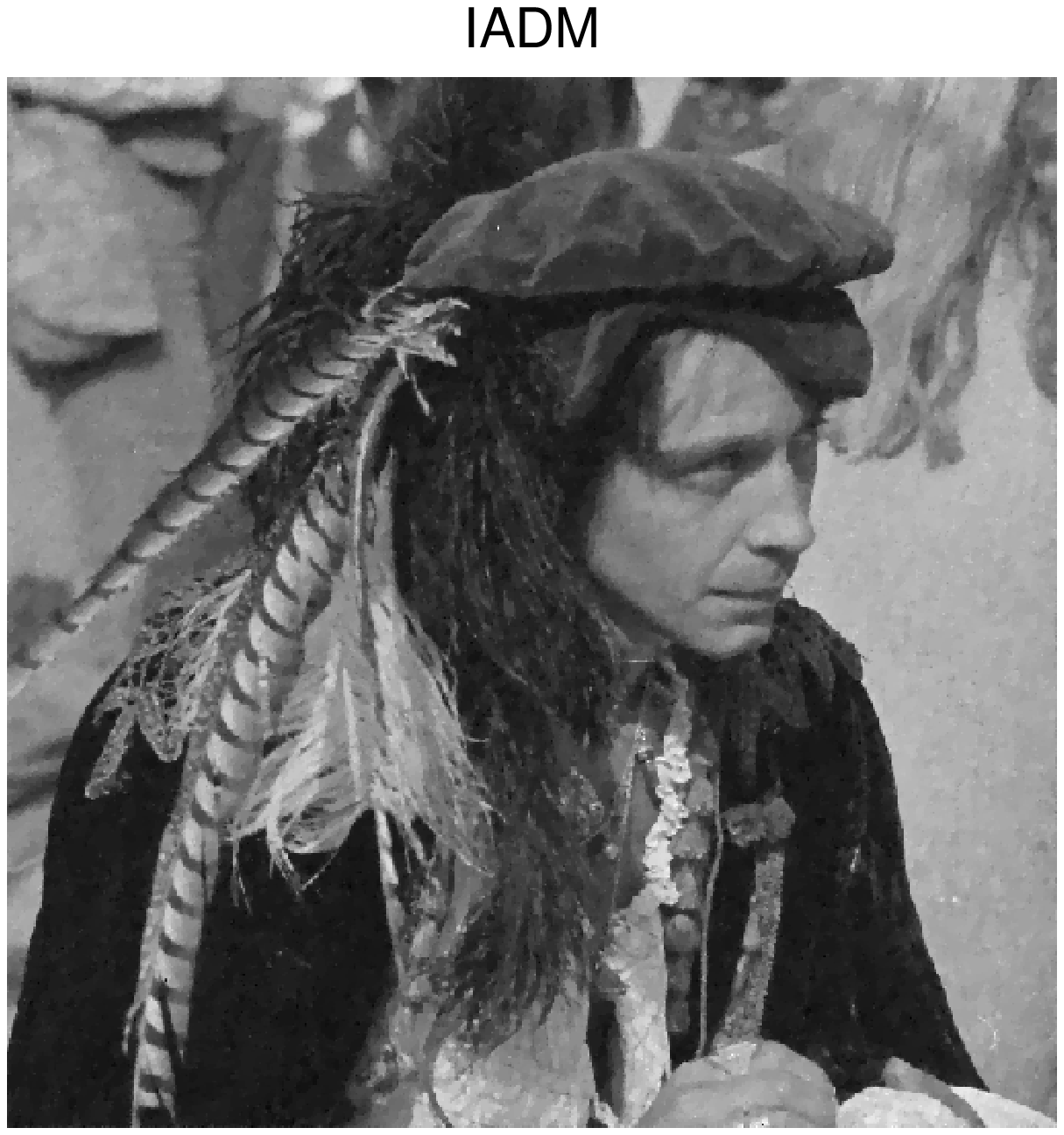}\hspace{.2cm}
\caption{{\footnotesize Original and recovered images by IADM and
TwIST. From top to bottom: brain 1, brain 2, cameraman, lena and man.  }}\label{figure4}
\end{figure}

\begin{figure}[htbp]
\vspace{-0.0cm}  \centering \vspace{0.5cm}
\includegraphics[scale = .285]{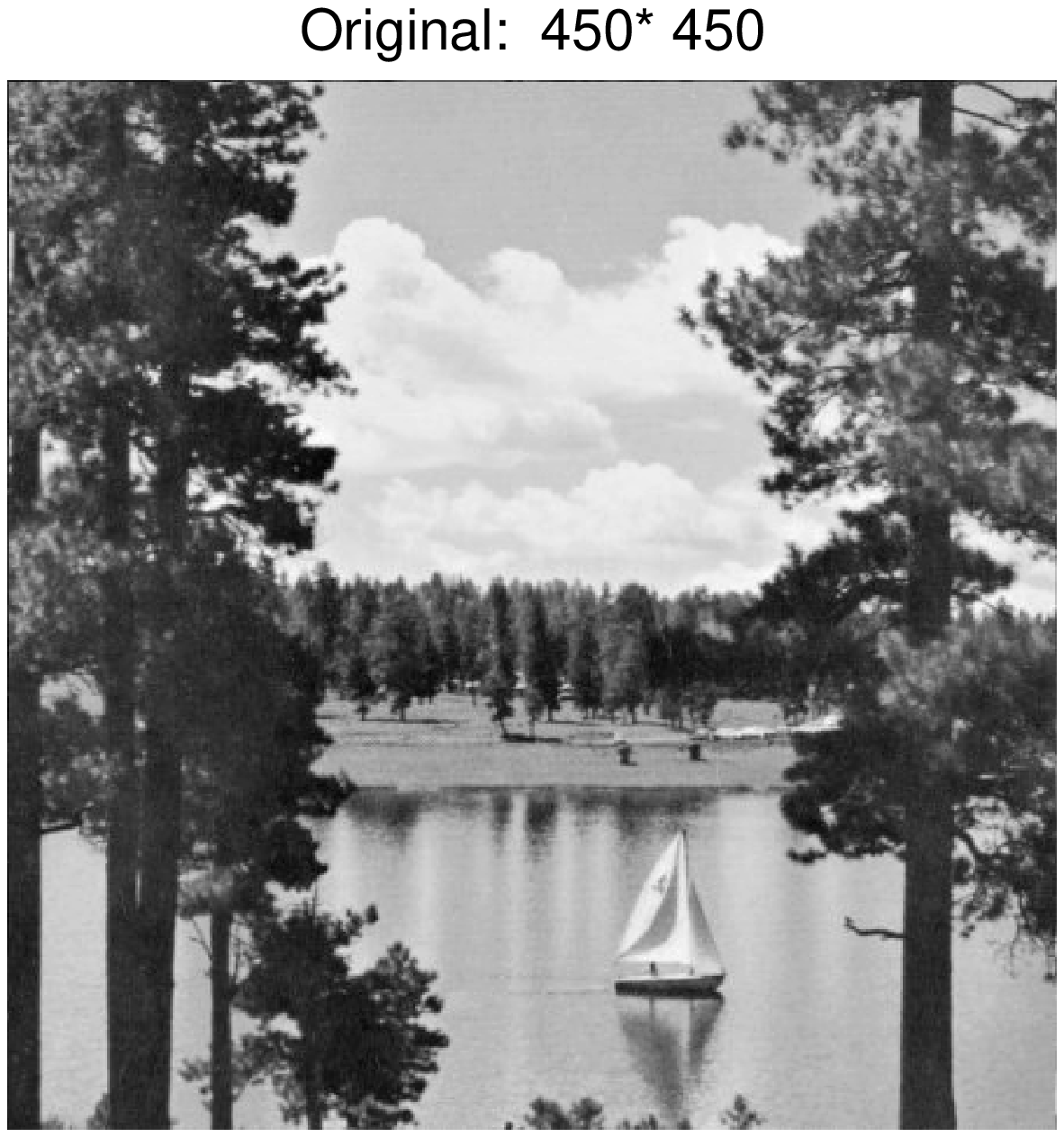}\hspace{.2cm}
\includegraphics[scale = .285]{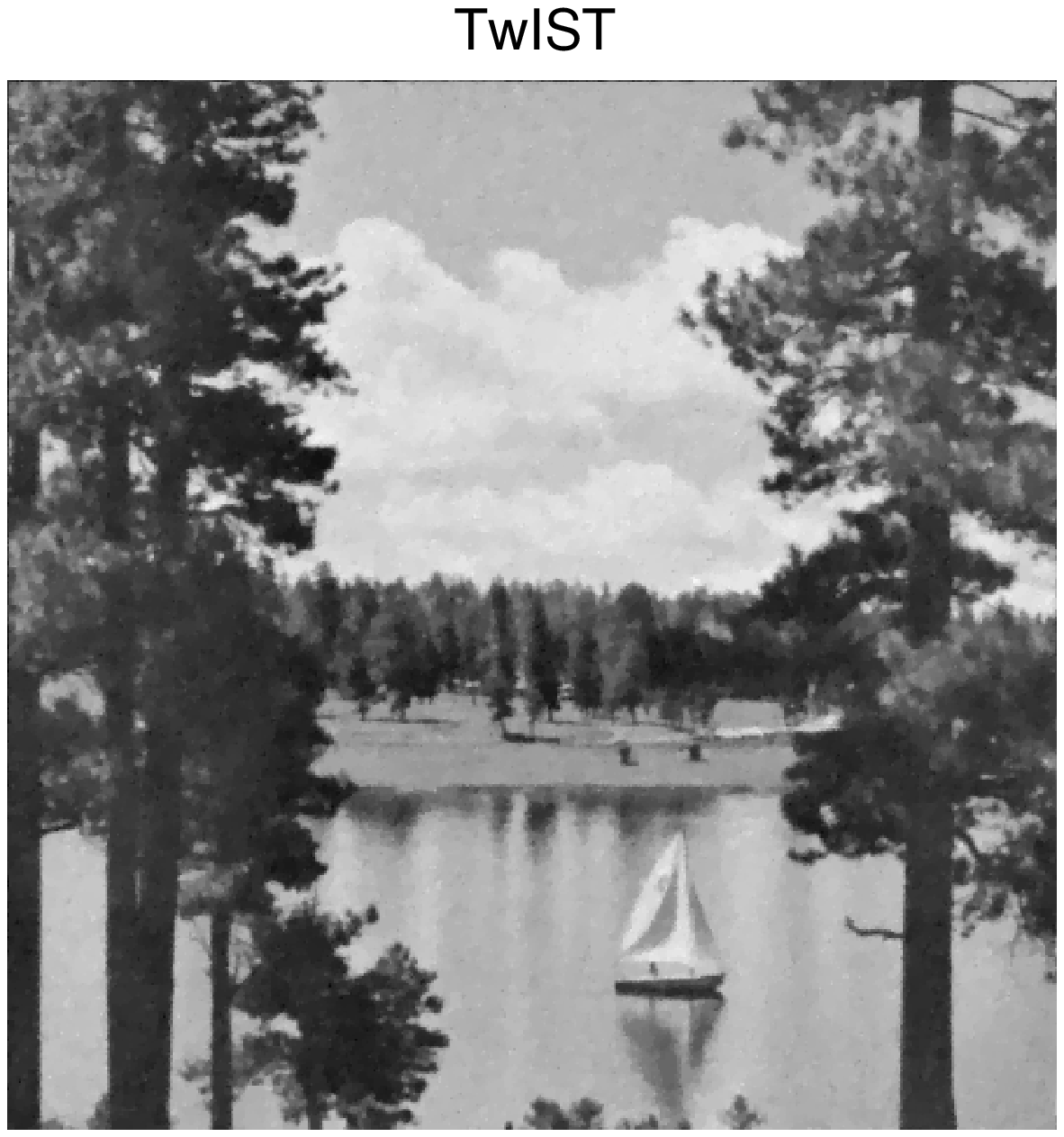}\hspace{.2cm}
\includegraphics[scale = .285]{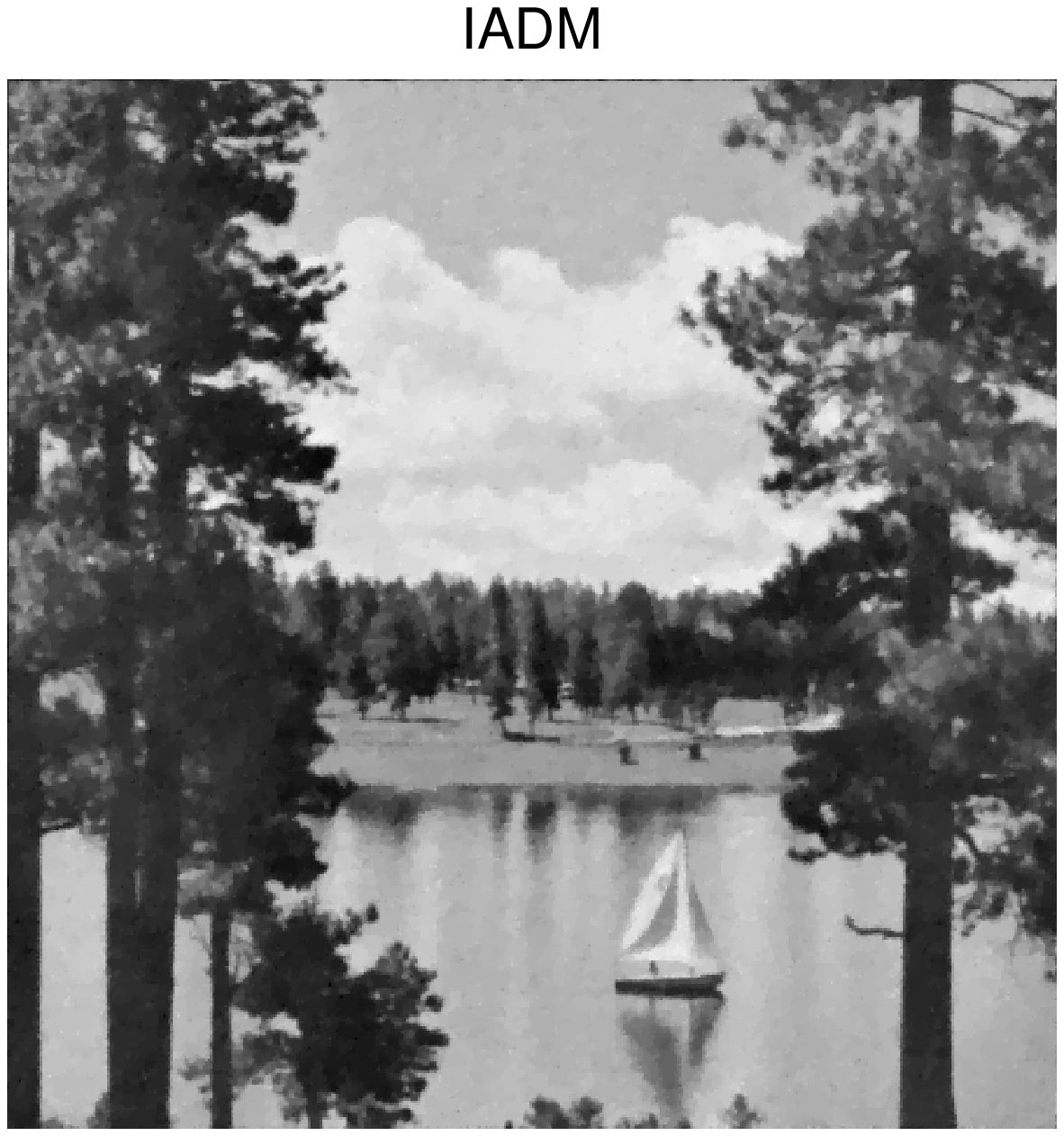}\hspace{.2cm}\\
\vspace{0.5cm}
\includegraphics[scale = .28]{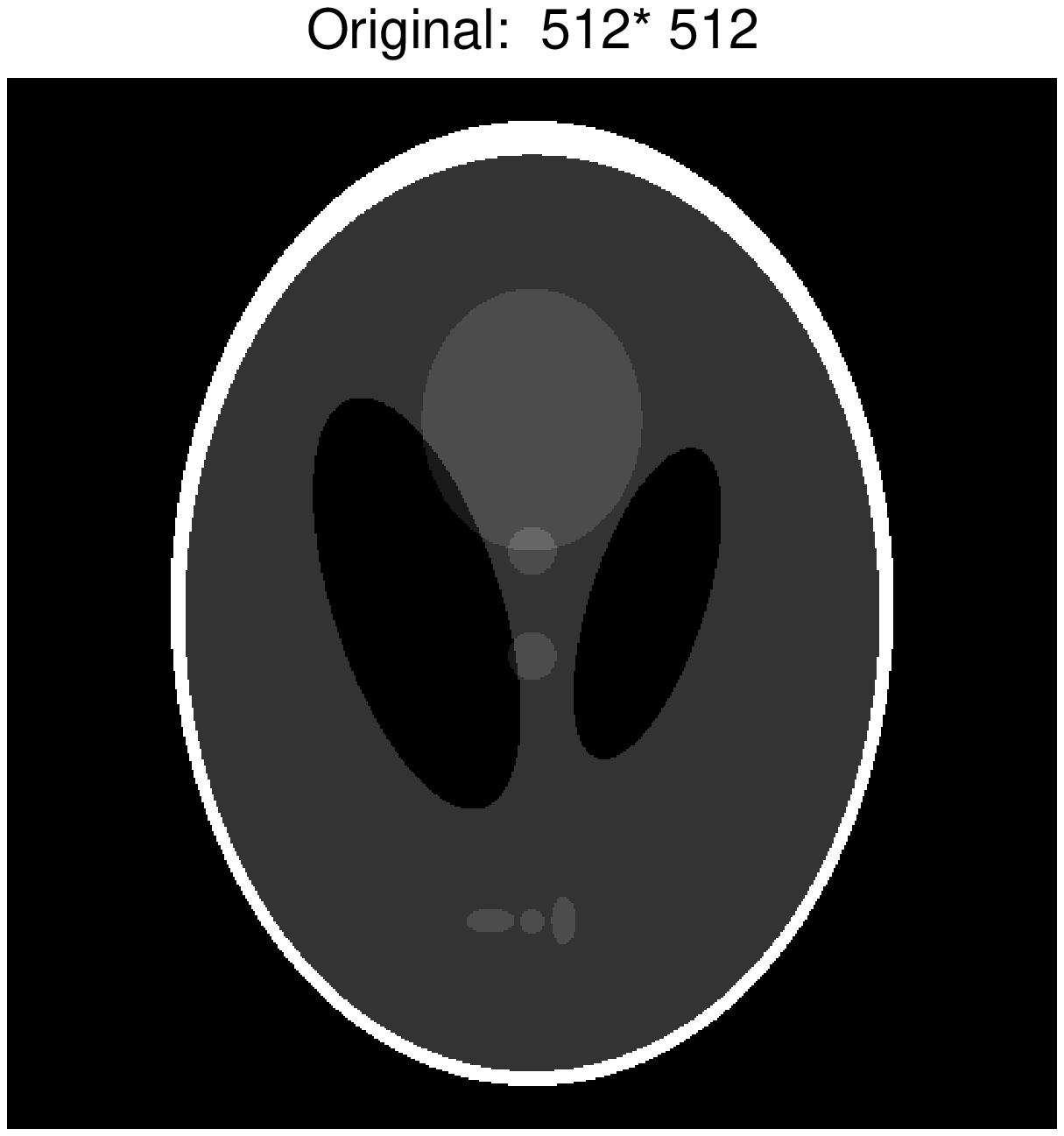}\hspace{.2cm}
\includegraphics[scale = .28]{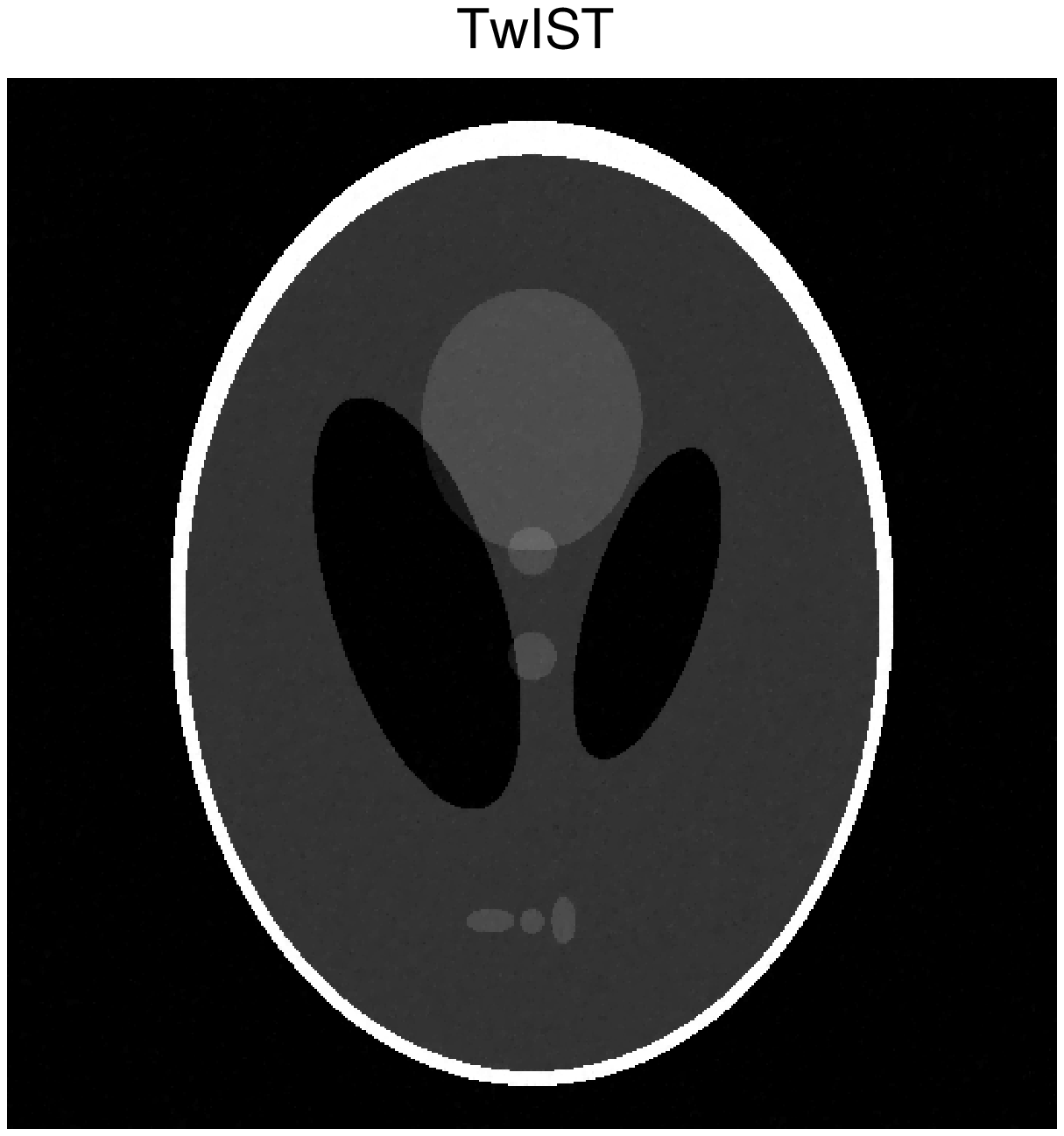}\hspace{.2cm}
\includegraphics[scale = .28]{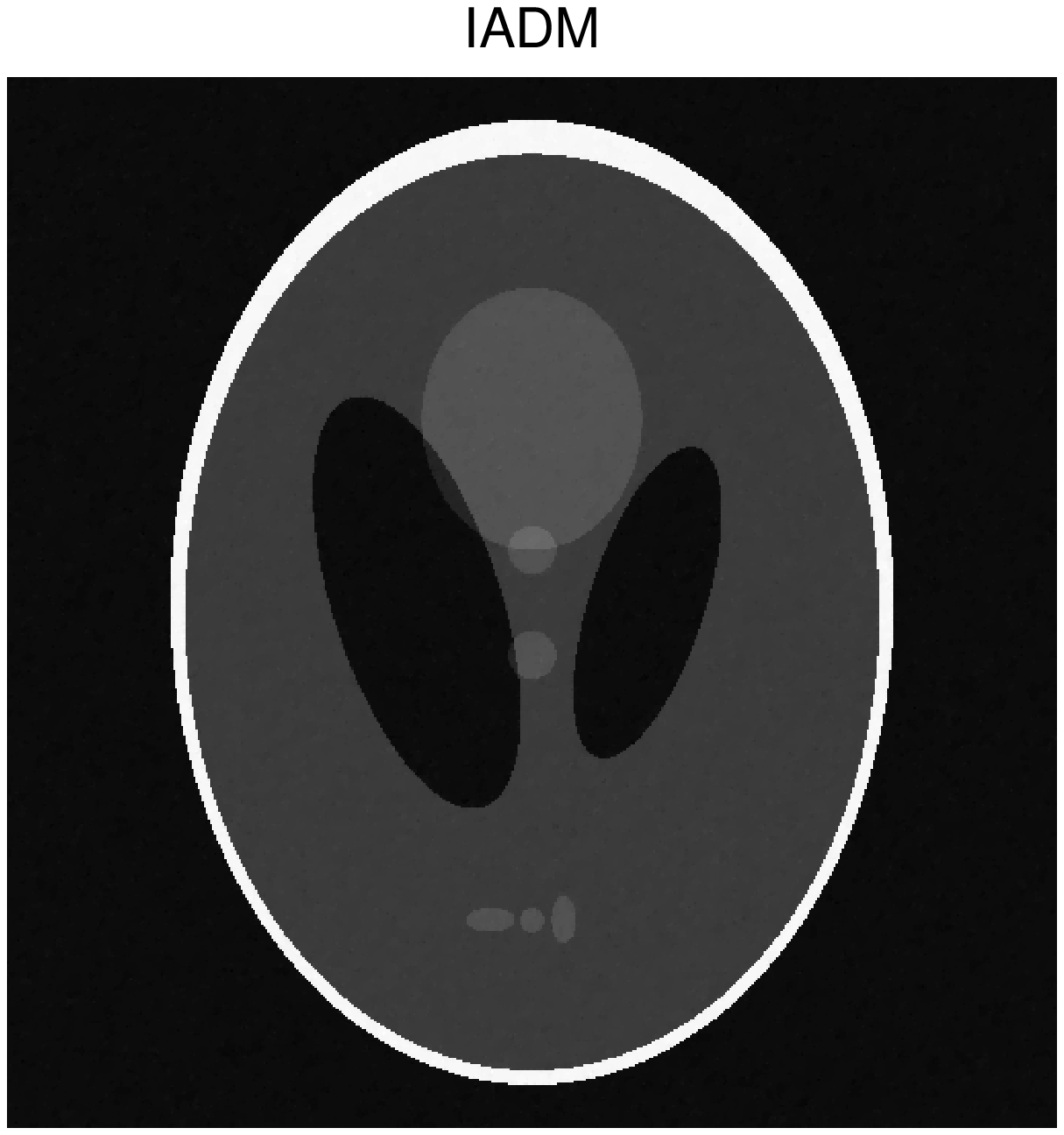}\hspace{.2cm}\\
\vspace{0.5cm}
\includegraphics[scale = .28]{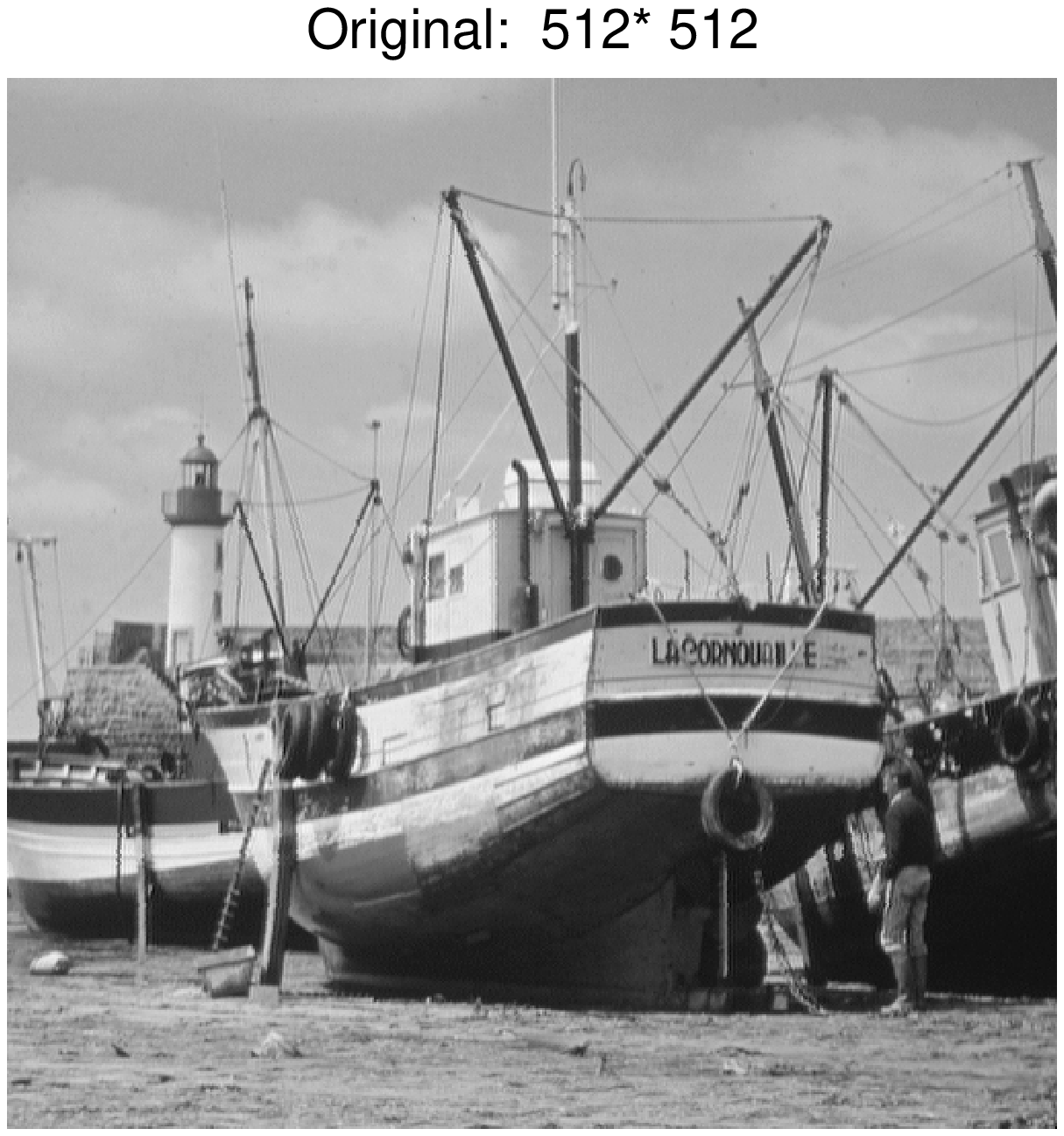}\hspace{.2cm}
\includegraphics[scale = .28]{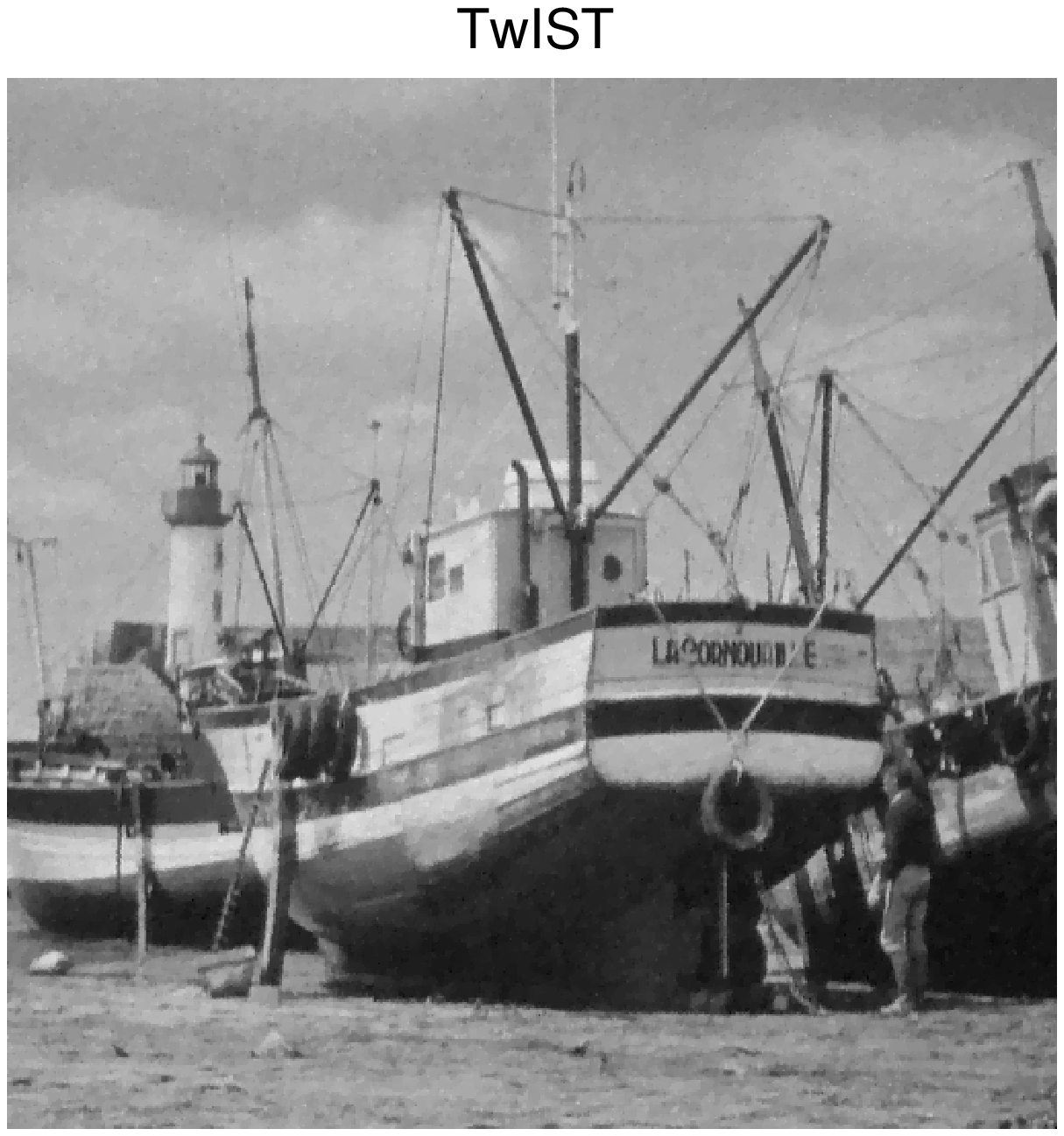}\hspace{.2cm}
\includegraphics[scale = .28]{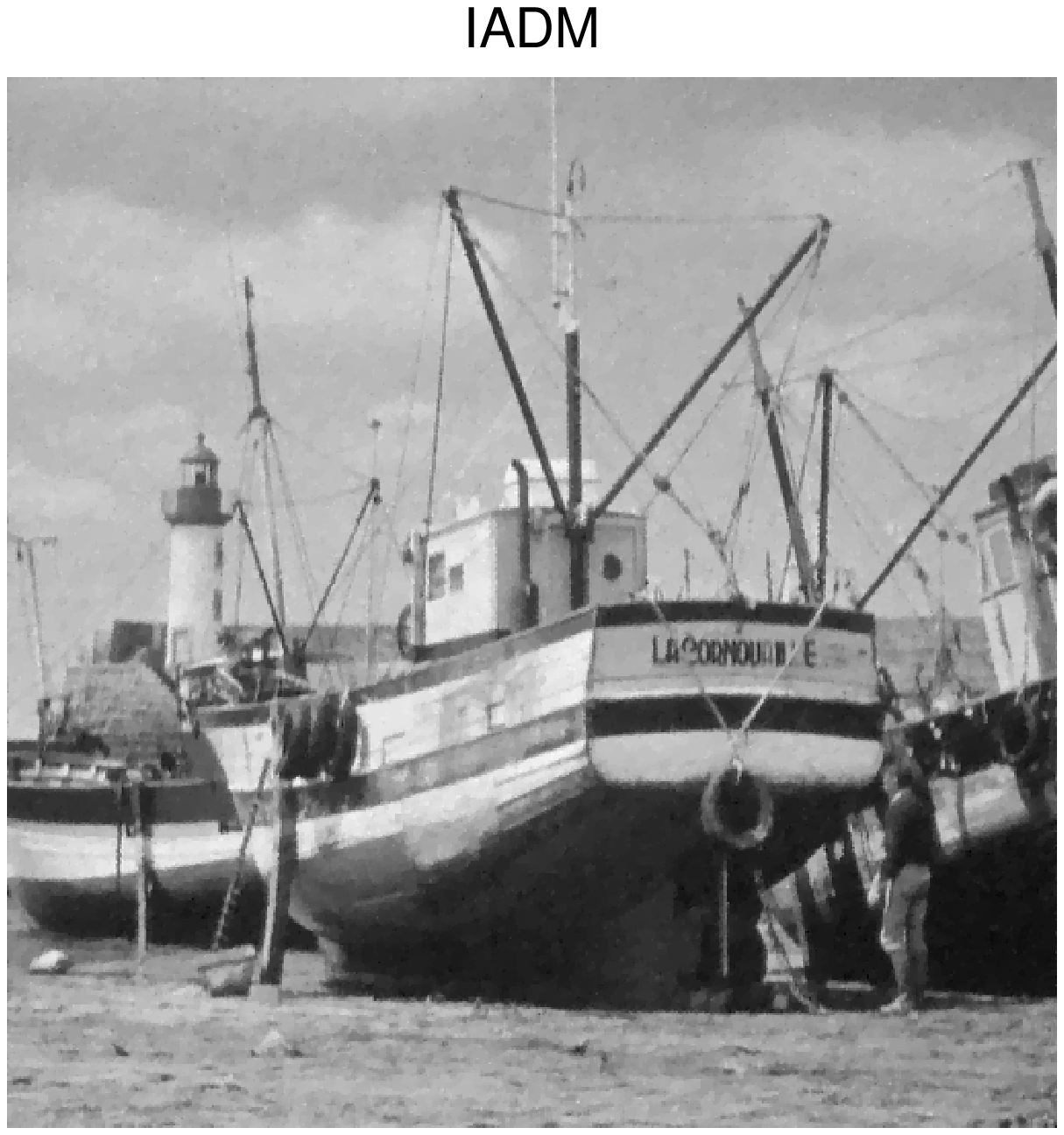}\hspace{.2cm}
\caption{{\footnotesize Original and recovered images by IADM and
TwIST. From top to bottom: sailboat, sheppon and boat.}}\label{figure5}
\end{figure}

Besides the Shepp-Logan phantom image, we also tested Cameraman, Lena,
Boat, Sailboat, as well as two brain images.  In this experiment, we simply set
$\mu=500$ and keep all other parameters unchanged.   The
original and the recovered images by TwIST and IADM are
given in Figures \ref{figure4} and \ref{figure5}, and detailed
results including relative errors (RE), CPU time (Time), final objective
function values (Obj), and the number of iterations (Iter) are presented in
Table \ref{table2}. It can be seen from Table \ref{table2} that
IADM attained comparable or better image quality in less CPU seconds.
For each test, IADM consumed more iterations while
the CPU time is less because the per-iteration cost of IADM is
much less than that of TwIST. Specifically, the per-iteration cost
of IADM contains two matrix-vector multiplications and two FFTs, while
TwIST needs to solve a TV denoising problem at each iteration.
In addition, IADM always attained smaller function values.
In summary, the comparison results indicate that
IAMD performs favorably and can be competitive with
the state-of-the-art algorithm TwIST.

\begin{table}
\begin{center}
\scriptsize{ \caption{Comparison results of IADM and TwIST with different images.}
\begin{tabular}{l|l|rrrr|rrrr}
\hline\hline \multicolumn{2}{c|}{ } & \multicolumn{4}{c|}{TwIST } &
\multicolumn{4}{c}{ IADM} \\ \hline
Images &Size & Iter  &  RE  &  Time  &  Obj  & Iter  &  Re  &  Time  & Obj\\
\hline
brain 1      &$128\times 128$ &     52 &  14.01\% & 34.25s & 4.7831e+002 &    208 & 13.64\% &  20.41s & 4.5478e+002\\
brain 2      &$256\times 256$ &     48 &  9.59\% & 90.22s & 1.6397e+003 &    176 & 9.45\% &  61.67s &  1.5665e+003\\
cameraman &$256\times 256$ &     56 &  5.71\% & 122.22s & 2.9068e+003 &    257 & 5.59\% &  118.67s & 2.7822e+003\\
lena      &$256\times 256$ &     53 &  4.93\% & 121.06s & 2.3656e+003 &    205 & 5.01\% &  92.11s &  2.2627e+003\\
man       &$512\times 512$ &     59 &  8.54\% & 423.38s & 1.0617e+004 &    262 & 8.57\% &  400.81s & 1.0122e+004 \\
sailboat  &$450\times 450$ &     57 &  4.91\% & 361.58s & 7.9220e+003 &   245& 4.98\% &  260.73s & 7.5960e+003\\
sheppon   &$512\times 512$ &     42 &  2.62\% & 335.36s & 4.4496e+003 &   135& 2.09\% &  217.08s & 4.2317e+003\\
boat      &$512\times 512$ &     53 &  4.37\% & 477.17s & 8.7306e+003 &   200& 4.34\% &  384.61s & 8.3312e+003\\
barbara   &$512\times 512$ &     56 &  9.83\% & 493.20s & 1.3469e+004 &   292& 9.80\% &  550.39s & 1.2814e+004\\
\hline\hline
\end{tabular} \label{table2}
}
\end{center}
\end{table}

%

\section{Concluding Remarks}\label{ConcludingRemarks}

In this paper, we proposed a Fast alternating minimization algorithm for
Total Variation image reconstruction from  Compressive Sensing data (FTVCS).
The per-iteration cost of FTVCS includes a linear time shrinkage
operation, two matrix-vector multiplications and two FFTs.
To overcome the difficulty caused by large penalty parameter in FTVCS,
we have also developed an Inexact Alternating Direction Method (IADM)
based on  linearization and proximal techniques. Our experimental results
indicate that IADM indeed performs better than FTVCS and is comparable with
the state-of-the-art algorithm TwIST for solving TV
reconstruction models.

Given the promising performance of IADM and the wide applications of TV models, we believe that
it is worthwhile  to further accelerate
IADM via certain line search strategy.  In both FTVCS and IADM, we used a linearization technique and FFTs to obtain a new point.
A possible improvement is to solve the $u$-subproblem of \eqref{ALP} by using certain gradient methods,
e.g., gradient descent method with BB steplengths \cite{BB} and non-monotone line search.
This should be interesting for further investigations.




\appendix

\section{Proof of Theorems \ref{conver}, \ref{linear} and \ref{qline}}\label{AppendixA}


The purpose of this appendix is to establish   convergence
properties of Algorithm \ref{ftvdcsalg} for a fixed $\beta>0$. For convenience, we define some notation.
For fixed $\beta>0$, the 2-dimensional (2D)
shrinkage operator $s:\R^2\rightarrow\R^2$ is defined as
\begin{equation}\label{sop}
s(\alpha)\triangleq\alpha-\mathcal{P_{\mathcal{B}}}(\alpha)=\max\left\{\|\alpha\|-1/\beta,0\right\}\frac{\alpha}{\|\alpha\|},
\end{equation}
where $\mathcal{P_{\mathcal{B}}}(\cdot):\R^2\rightarrow\R^2$
is the projection onto the closed disc
$\mathcal{B}\triangleq\{\alpha\in\R^2:\|\alpha\|\leq 1/\beta\}$, and the convention $0\cdot(0/0)=0$ is followed.
For vectors $u,v\in \R^N$, $N\geq 1$, we define
$S(u;v):\R^{2N}\rightarrow\R^{2N}$ by
$$
S(u;v)=(s(\alpha_1);\ldots;s(\alpha_N)), \ \ \mathrm{where} \
\alpha_i=\left(
\begin{array}{c}
u_i\\
v_i
\end{array}
\right),
$$
i.e., $S$ applies 2D shrinkage to each pair $(u_i;v_i)$, for $i=1,2,\ldots,N$.
From the definition of $s(\cdot)$, it is easy to see that
(\ref{fixu}) or (\ref{shrink}) can be rewritten as
$\mathbf{w}_i=s(D_iu)$.
The following result shows that the operator $s$ is non-expansive.
\begin{lemma}[\cite{YANG1}]\label{nonexp}
For any $a, b\in \R^2$, it holds that
$$
\|s(a)-s(b)\|^2\leq
\|a-b\|^2-\|\mathcal{P}_{\mathcal{B}}(a)-\mathcal{P}_\mathcal{B}(b)\|^2,
$$
Furthermore, if $\|s(a)-s(b)\| =\|a-b\|$, then
$s(a)-s(b)=a-b$.
\end{lemma}

Since the objective function in (\ref{atv}) is convex, bounded
below, and coercive (i.e., its goes to infinity as
$\|(w,u)\|\rightarrow \infty)$, it has at least one
minimizer $(w^*,u^*)$ that cannot be decreased by the
alternating minimization scheme \eqref{shrink}-\eqref{foru} and thus must satisfy
\begin{equation}\label{solution}
\left \{
\begin{array}{l}
w^*=S(D^{(1)}u^*;D^{(2)}u^*) \ (\triangleq S(Du^*)), \\
(D^\top  D+\frac{\mu}{\beta}A^\top  A)u^*=D^\top  w^*+\frac{\mu}{\beta}A^\top  f.
\end{array}
\right.
\end{equation}
By using the shrinkage operator, we can rewrite the iteration of Algorithm \ref{ftvdcsalg}
as
\begin{equation}\label{uk1}
\left \{
\begin{array}{l}
w^{k+1}=S(D^{(1)}u^k;D^{(2)}u^k) \ (\triangleq S(Du^k)), \\
(D^\top  D+\frac{\mu}{\beta\tau}I)u^{k+1}=D^\top  w^{k+1}+\frac{\mu}{\beta\tau}(u^k-\tau
g_k).
\end{array}
\right.
\end{equation}

In the following, we show that (\ref{uk1}) converges to \eqref{solution}.
The following matrices will be used in our analysis:
$$
M=D^\top  D+\frac{\mu}{\beta}A^\top  A, \ \ H=D^\top  D+\frac{\mu}{\beta\tau}I, \ \
\mbox{and} \ \ T=I-\tau A^\top  A.
$$
Assumption \ref{assum1} ensures  the non-singularity of $M$, while
$H^{-1}$ is always well defined under the circumstance. Simple manipulation shows that
$H-M=\eta^2T$, where $\eta\triangleq\sqrt{\frac{\mu}{\beta\tau}}$.
With these definitions, (\ref{solution}) and (\ref{uk1}) can be, respectively,  simplified as
$$
\left \{
\begin{array}{l}
w^*= S(D  u^*), \\
v^*=\eta T u^*+\eta\tau A^\top  f,\\
Hu^*=D^\top  w^*+\eta v^*,
\end{array}
\right.
\quad \mathrm{and}\quad
\left \{
\begin{array}{l}
w^{k+1}= S(D  u^k), \\
v^{k+1}=\eta T u^k+\eta\tau A^\top  f,\\
Hu^{k+1}=D^\top  w^{k+1}+\eta v^{k+1}.
\end{array}
\right.
$$
To further simplify the above equations, we define
$$
h(w;v)=DH^{-1}\left(
\begin{array}{c}
D\\
\eta I
\end{array}
\right)^\top   \left(
\begin{array}{c}
w\\
v
\end{array}
\right)
\quad \mathrm{and}\quad
p(w;v)=\eta TH^{-1}\left(
\begin{array}{c}
D\\
\eta I
\end{array}
\right)^\top   \left(
\begin{array}{c}
w\\
v
\end{array}
\right).
$$
Hence, the solution and iteration systems can be, respectively, rewritten
as
\begin{equation}\label{solutionrew11}
\left \{
\begin{array}{l}
w^*= S\circ h(w^*;v^*), \\
v^*= p(w^*;v^*)+\eta\tau A^\top  f,\\
Hu^*=D^\top  w^*+\eta v^*,
\end{array}
\right.
\end{equation}
and
\begin{equation}\label{uk1rew22}
\left \{
\begin{array}{l}
w^{k+1}=S\circ h(w^k;v^k), \\
v^{k+1}=p(w^k;v^k)+\eta\tau A^\top  f,\\
Hu^{k+1}=D^\top  w^{k+1}+\eta v^{k+1},
\end{array}
\right.
\end{equation}
where   ``$\circ$" denotes operator composition.
Furthermore, we define
$$
q(w;v)=\left(
\begin{array}{c}
S\circ h(w;v)\\
p(w;v)
\end{array}
\right) + \left(
\begin{array}{c}
0\\
\eta\tau A^\top  f
\end{array}
\right).
$$
Then (\ref{solutionrew11}) and (\ref{uk1rew22}) become
$$
\left \{
\begin{array}{l}
(w^*;v^*)=q(w^*;v^*)\\
Hu^*=D^\top  w^*+\eta v^*,
\end{array}
\right.
\quad \mathrm{and}\quad
\left \{
\begin{array}{l}
(w^{k+1};v^{k+1})=q(w^k;v^k)\\
Hu^{k+1}=D^\top  w^{k+1}+\eta v^{k+1},
\end{array}
\right.
$$

\begin{lemma}\label{lem41}
$q(w;v)$ is non-expansive.
\end{lemma}
\begin{proof}
Given $(w^1;v^1)$ and $(w^2;v^2)$, it holds that
\begin{eqnarray*}
\|q(w^1;v^1)-q(w^2;v^2)\|^2 &=& \left \|
\begin{array}{c}
S\circ h(w^1;v^1)-S\circ h(w^2;v^2)\\
p(w^1;v^1)-p(w^2;v^2)
\end{array}
\right \|^2\\
&\leq&\|h(w^1;v^1)-h(w^2;v^2)\|^2+\|p(w^1;v^1)-p(w^2;v^2)\|^2\\
&=& \left \|R \left(
\begin{array}{c}
w^1-w^2\\
v^1-v^2
\end{array}
\right )\right \|^2,
\end{eqnarray*}
where ``$\leq$" comes from the non-expansive of $s(\cdot)$ and
$$
R\triangleq\left (
\begin{array}{c}
D\\
\eta T
\end{array}
\right )H^{-1} \left (
\begin{array}{c}
D\\
\eta I
\end{array}
\right )^\top  .
$$
It is easy to verify that
\begin{eqnarray*}
R^\top  R &=&\left (
\begin{array}{c}
D\\
\eta I
\end{array}
\right )H^{-1}(D^\top  D+\eta^2T^2)H^{-1}\left (
\begin{array}{c}
D\\
\eta I
\end{array}
\right )^\top  \\
&=& \left (
\begin{array}{c}
D\\
\eta I
\end{array}
\right
)H^{-1}\left(H-\frac{\mu}{\beta}(2A^\top  A-\tau(A^\top  A)^2)\right)H^{-1}\left
(
\begin{array}{c}
D\\
\eta I
\end{array}
\right )^\top  \\
&=& \left (
\begin{array}{c}
D\\
\eta I
\end{array}
\right)H^{-1}\left (
\begin{array}{c}
D\\
\eta I
\end{array}
\right)^\top  -\frac{\mu}{\beta} \left (
\begin{array}{c}
D\\
\eta I
\end{array}
\right)H^{-1}(2A^\top  A-\tau(A^\top  A)^2)H^{-1}\left (
\begin{array}{c}
D\\
\eta I
\end{array}
\right )^\top  .
\end{eqnarray*}
Recall that we require $0<\tau<2/\lambda_{\max}(A^\top A)$, which ensures the positive
semi-definiteness of $A^\top  A-\tau(A^\top  A)^2$. Therefore,
\begin{equation}\label{nonequ}
\|q(w^1;v^1)-q(w^2;v^2)\|^2\leq \left(
\begin{array}{c}
w^1-w^2\\
v^1-v^2
\end{array}
\right )^\top  \left (
\begin{array}{c}
D\\
\eta I
\end{array}
\right )H^{-1}\left (
\begin{array}{c}
D\\
\eta I
\end{array}
\right )^\top  \left(
\begin{array}{c}
w^1-w^2\\
v^1-v^2
\end{array}
\right )\leq \left\|
\begin{array}{c}
w^1-w^2\\
v^1-v^2
\end{array}
\right \|^2,
\end{equation}
which shows that $q(w;v)$ is non-expansive.
\end{proof}

\begin{lemma}\label{lem42}
Equality holds in (\ref{nonequ}) if and only if
$$
q(w^1;v^1)-q(w^2;v^2)=\left(
\begin{array}{c}
w^1-w^2\\
v^1-v^2
\end{array}
\right ).
$$
\end{lemma}
\begin{proof}
We note that in the proof of Lemma \ref{lem41} there exist three
``$\leq$". Thus, equality holds in \eqref{nonequ} only when all the three inequalities become
``=". For simplicity, we let $dw=w^1-w^2$ and $dv=v^1-v^2$.
\begin{enumerate}
  \item The first ``$\leq$" becomes ``$=$" if and only if
$$
S\circ h(w^1;v^1)-S\circ h(w^2;v^2)=h(w^1;v^1)-
h(w^2;v^2)=DH^{-1}\left (
\begin{array}{c}
D\\
\eta I
\end{array}
\right )^\top  \left (
\begin{array}{c}
dw\\
dv
\end{array}
\right ).
$$
  \item The second ``$\leq$" becomes ``$=$" if and only if
$$
\left (
\begin{array}{c}
dw\\
dv
\end{array}
\right )^\top   \left (
\begin{array}{c}
D\\
\eta I
\end{array}
\right )H^{-1}(2A^\top  A-\tau(A^\top  A)^2)H^{-1}\left (
\begin{array}{c}
D\\
\eta I
\end{array}
\right )^\top  \left (
\begin{array}{c}
dw\\
dv
\end{array}
\right )=0.
$$
  \item Let $U$ be orthonormal and
$$
\left (
\begin{array}{c}
D\\
\eta I
\end{array}
\right )^\top  H^{-1}\left (
\begin{array}{c}
D\\
\eta I
\end{array}
\right ) = U^\top  \Lambda U
$$
be its eigenvalue decomposition. The third ``$\leq$" becomes ``$=$" if
and only if
$$
\sum_{i}\lambda_i\left ( U \left (
\begin{array}{c}
dw\\
dv
\end{array}
\right )\right )_i^2=\sum_{i}\left ( U \left (
\begin{array}{c}
dw\\
dv
\end{array}
\right )\right )_i^2.
$$
Since $0\leq \lambda_i\leq 1$, the above equality holds only when
$$
\lambda_i\left ( U \left (
\begin{array}{c}
dw\\
dv
\end{array}
\right )\right )_i^2=\left ( U \left (
\begin{array}{c}
dw\\
dv
\end{array}
\right )\right )_i^2, \ \ \ \forall \ i.
$$
Therefore,
$$
\Lambda U \left (
\begin{array}{c}
dw\\
dv
\end{array}
\right )=U \left (
\begin{array}{c}
dw\\
dv
\end{array}
\right ),
$$
and thus
\begin{equation}\label{eq41}
\left (
\begin{array}{c}
D\\
\eta I
\end{array}
\right )H^{-1}\left (
\begin{array}{c}
D\\
\eta I
\end{array}
\right )^\top   \left (
\begin{array}{c}
dw\\
dv
\end{array}
\right )=\left (
\begin{array}{c}
dw\\
dv
\end{array}
\right ).
\end{equation}
\end{enumerate}
From 1 and (\ref{eq41}), we have
\begin{equation}\label{eq42}
S\circ h(w^1;v^1)-S\circ h(w^2;v^2)=DH^{-1} \left (
\begin{array}{c}
D\\
\eta I
\end{array}
\right )^\top   \left (
\begin{array}{c}
dw\\
dv
\end{array}
\right )=dw.
\end{equation}
From (\ref{eq41}), the equality in 2 is equivalent to
$$
dv^\top  (2A^\top  A-\tau(A^\top  A)^2)dv=0.
$$
Let $U^\top  \Lambda U=A^\top  A$ be the eigenvalue decomposition of $A^\top  A$.
The above equation is equivalent to
$$
dv^\top  U^\top  (2\Lambda-\tau \Lambda^2)Udv=0 \quad\mathrm{or }\quad
\sum_i(2\lambda_i-\tau\lambda_i^2)(Udv)_i^2=0.
$$
Since $2\lambda_i-\tau\lambda_i^2\geq 0$, we have
$
(2\lambda_i-\tau\lambda_i^2)(Udv)_i^2=0,    \forall \ i.
$
If $\lambda_i\neq 0$, then from the choice of $\tau$ we have
$2\lambda_i-\tau\lambda_i^2>0$, and  thus $(Udv)_i=0$. Therefore,
$\Lambda U dv=0$ and $
A^\top  Adv=U^\top  \Lambda U dv=0.
$
Sum the above discussions up, we have
\begin{eqnarray*}
q(w^1;v^1)-q(w^2;v^2) &=& \left (
\begin{array}{c}
S\circ h(w^1;v^1)-S\circ h(w^2;v^2)\\
p(w^1;v^1)-p(w^2;v^2)
\end{array}
\right )\\
 &=& \left (
\begin{array}{c}
dw\\
T\cdot dv
\end{array}
\right ) = \left (
\begin{array}{c}
dw\\
dv-A^\top  A dv
\end{array}
\right ) = \left(
\begin{array}{c}
w^1-w^2\\
v^1-v^2
\end{array}
\right ),
\end{eqnarray*}
where the first equality is from the definition of $q(\cdot;\cdot)$; the second one is
from (\ref{eq42}), the definition of $p$ and (\ref{eq41}); the third
one is from the definition of $T$; and the final one is from
$A^\top  dv=0$. This completes the proof.
\end{proof}

\begin{corollary}
Suppose $(w^*;v^*)$ is a fixed point of $q$, i.e.,
$(w^*;v^*)=q(w^*;v^*)$. Then for any $(w;v)$ it holds
$$
\|q(w;v)-q(w^*;v^*)\|<\|(w;v)-(w^*;v^*)\|
$$
unless $(w;v)$ is also a fixed point of $q(\cdot;\cdot)$.
\end{corollary}

Based on the above lemmas, now we are ready to give the proofs of Theorems
\ref{conver}, \ref{linear} and \ref{qline}.

\begin{proof}
(Theorem \ref{linear})
First, the convergence of
$(w^k,u^k)$ to $(w^*,u^*)$ can be established using exactly the same arguments as in
Theorem 3.4 in \cite{YANG1}. The convergence of $u^k$ to $u^*$
follows from the convergence of $w^k$ to $w^*$ and $v^k$ to
$v^*$. Therefore, we omit the details.
\end{proof}

For any $i$, we let
$
h_i(w;v)=D_iH^{-1}(D^\top  w+\eta v): \R^{3n^2}\rightarrow
\R^2
$,
$
E=\{1,\ldots,n^2\}/L,
$,
where we recall that
$
L=\{i:\|D_iu^*\|=\|h_i(w^*,v^*)\|\leq 1/\beta\}
$,
and
\begin{equation}\label{eq51}
w=\min\{1/\beta-\|D_iu^*\|:i\in L\}>0.
\end{equation}

\begin{proof}
({\sc Theorem \ref{linear}}) From the non-expansive of $s(\cdot)$, for
each $i$, it holds
\begin{equation}\label{eq52}
\|\mathbf{w}_i^{k+1}-\mathbf{w}_i^*\|=\|s\circ h_i(w^k;v^k)-s\circ
h_i(w^*;v^*)\|\leq \| h_i(w^k;v^k)-  h_i(w^*;v^*)\|.
\end{equation}
Suppose that at iteration $k$ there exist at least one index $i\in
L$ such that $\mathbf{w}_i^{k+1}=s\circ h_i(w^k;v^k)\neq 0$. Then
$\| h_i(w^*;v^*)\|\leq 1/\beta$, $\| h_i(w^k;v^k)\|> 1/\beta$, and
$\mathbf{w}_i^{*}=s\circ h_i(w^*;v^*)= 0$. Therefore,
\begin{eqnarray}
\|\mathbf{w}_i^{k+1}-\mathbf{w}_i^*\| &=& \|s\circ
h_i(w^k;v^k)\|^2=(\| h_i(w^k;v^k)\|-1/\beta)^2\label{eq53}\\
&\leq&
\left[\|h_i(w^k;v^k)-h_i(w^*;v^*)\|-(1/\beta-\|h_i(w^*;v^*)\|\right]^2\nonumber\\
&\leq&
\|h_i(w^k;v^k)-h_i(w^*;v^*)\|^2-(1/\beta-\|h_i(w^*;v^*)\|)^2\nonumber\\
&\leq& \|h_i(w^k;v^k)-h_i(w^*;v^*)\|^2-\omega^2,\nonumber
\end{eqnarray}
where the first ``$\leq$" is the triangular inequality, the second
one follows from the fact that
$\|h_i(w^k;v^k)-h_i(w^*;v^*)\|\geq1/\beta-\|h_i(w^*;v^*)\|>0$, and
the last  one used the definition of $\omega$ in (\ref{eq51}).
Combining with (\ref{eq51}) and (\ref{eq53}), we obtain
\begin{eqnarray*}
\left \|
\begin{array}{c}
w^{k+1}-w^*\\
v^{k+1}-v^*
\end{array}
\right \|^2&=&\sum_i
\|\mathbf{w}_i^{k+1}-\mathbf{w}_i^*\|^2+\|v^{k+1}-v^*\|^2\\
&\leq& \sum_i \|h_i(w^k;v^k)-h_i(w^*;v^*)\|^2-\omega^2+\|v^{k+1}-v^*\|^2\\
&=&
\|h_i(w^k;v^k)-h_i(w^*;v^*)\|^2-\omega^2+\|p(w^k;v^k)-p(w^*;v^*)\|\\
&\leq & \left \|
\begin{array}{c}
w^{k}-w^*\\
v^{k}-v^*
\end{array}
\right \|^2-\omega^2,
\end{eqnarray*}
where the second ``$\leq$" comes from the
non-expansiveness of
$$
\phi(w;v) =\left(
\begin{array}{c}
h(w;v)\\
p(w;v)
\end{array}
\right ),
$$
which  can be easily derived. Therefore, the number of iterations $k$
with $\mathbf{w}_i^{k+1}\neq 0$ does not exceed
$$
\frac{1}{\omega^2} \left \|
\begin{array}{c}
w^0-w^*\\
v^0-v^*
\end{array}
\right \|^2.
$$
This completes the proof of Theorem \ref{linear}.
\end{proof}

\begin{proof}
(Theorem \ref{qline}) From the iteration formulae for $u$ and
$(w;v)$, there holds
$$
u^{k+1}-u^*=H^{-1} \left (
\begin{array}{c}
D\\
\eta I
\end{array}
\right )^\top   \left (
\begin{array}{c}
w^{k+1}-w^*\\
v^{k+1}-v^*
\end{array}
\right ),
$$
and
\begin{eqnarray*}
\left \|
\begin{array}{c}
w^{k+1}-w^*\\
v^{k+1}-v^*
\end{array}
\right \|^2 &=& \|q(w^k;v^k)-q(w^*;v^*)\|^2 \leq \left \|
\begin{array}{c}
D(u^k-u^*)\\
\eta T(u^k-u^*)
\end{array}
\right \|^2 = \left \| R \left (
\begin{array}{c}
w^{k}-w^*\\
v^{k}-v^*
\end{array}
\right )\right \|^2.
\end{eqnarray*}
Considering the finite convergence of $\mathbf{w}_i$, $i\in L$, we
have
$$
\left \|
\begin{array}{c}
w_E^{k+1}-w_E^*\\
v^{k+1}-v^*
\end{array}
\right \|^2 \leq \rho((R^\top  R)_{EE}) \left \|
\begin{array}{c}
w_E^{k}-w_E^*\\
v^{k}-v^*
\end{array}
\right \|^2,
$$
where $(R^\top  R)_{EE}$ is a sub-matrix of $R^\top  R\in
\R^{3n^2\times3n^2}$ formed by throwing away certain rows (with indexes $\cup_{i\in L}\{i,i+n^2\}$)
and corresponding columns. Multiplying
$(D;\eta T)$ to the recursion of $u^{k+1}-u^*$, we get
\begin{eqnarray*}
\|u^{k+1}-u^*\|^2_{D^\top  D+\eta^2T^2} &=& \left (
\begin{array}{c}
w^{k+1}-w^*\\
v^{k+1}-v^*
\end{array}
\right )^\top   R^\top  R\left (
\begin{array}{c}
w^{k+1}-w^*\\
v^{k+1}-v^*
\end{array}
\right )\\
&=& \rho((R^\top  R)_{EE}) \left \|
\begin{array}{c}
w^{k+1}-w^*\\
v^{k+1}-v^*
\end{array}
\right \|^2 \leq \rho((R^\top  R)_{EE})\|u^k-u^*\|_{D^\top  D+\eta^2T^2}^2,
\end{eqnarray*}
which shows that $\{u^k\}$ converges q-linearly.
\end{proof}

\section{Convergence of Algorithm \ref{iadmalg}}\label{AppendixB}
In this section, we clarify the relationship between Algorithm \ref{iadmalg} and the
proximal  ADM approach proposed in \cite{BSHE}. The convergence of Algorithm \ref{iadmalg} follows
directly.

We briefly review the proximal ADM approach in \cite{BSHE} for structured  variational inequality (SVI) problems.
Let $M$ and $N$ be, respectively, $l\times n$ and $l\times m$ matrixes,
$\mathcal{X}\subset\R^n$ and $\mathcal{Y}\subset\R^m$ be nonempty closed convex sets,
and $f,g$ be given monotone operators. The SVI problem is to find
$u^*\in\Omega$ such that
\begin{equation}\label{vip}
(u-u^*)^\top  F(u^*)\geq 0, \ \ \forall \ u\in \Omega,
\end{equation}
where   $\Omega=\{(x,y):x\in\mathcal{X},y\in\mathcal{Y}, Mx+Ny=0\}$,
$$
u=\left(
\begin{array}{c}
x\\
y
\end{array}
\right ), \quad \text{and} \quad F(u)=\left(
\begin{array}{c}
f(x)\\
g(y)
\end{array}
\right ).
$$
Given $(x^k,y^k,\lambda^k)$, the proximal   ADM proposed in \cite{BSHE} iterates as follows
\begin{enumerate}
  \item Compute $x^{k+1}\in\mathcal{X}$ via solving
\begin{equation}\label{sub1}
(x'-x)^\top
\left\{f(x)-M^\top  [\lambda^k-h(Mx+Ny^k)]+R_k(x-x^k)\right\}\geq 0, \ \
\forall \ x'\in \mathcal{X}.
\end{equation}
  \item Compute $y^{k+1}\in\mathcal{Y}$ via  solving
\begin{equation}\label{sub2}
(y'-y)^\top
\left\{g(y)-N^\top  [\lambda^k-h(Mx^{k+1}+Ny)]+S_k(y-y^k)\right\}\geq 0,
\ \ \forall \ y'\in \mathcal{Y}.
\end{equation}
  \item Update $\lambda^{k+1}$ via
$$
\lambda^{k+1}=\lambda^k-h(Mx^{k+1}+Ny^{k+1}),
$$
\end{enumerate}
where $h>0$ is a parameter, $R_k$ and $S_k$ are symmetric positive semidefinite matrices.
Under mild assumptions, global convergence of this proximal ADM approach was established in \cite{BSHE}.
Simple manipulation shows that Algorithm \ref{iadmalg} is a special case of the proximal   ADM
approach described above by setting $R_k\equiv0$ in (\ref{sub1}), i.e., the the $w$-subproblems are solved exactly in Algorithm \ref{iadmalg}, and $S_k= \frac{1}{\tau}I-A^\top  A $ in \eqref{sub2}. Hence, the global convergence of   Algorithm \ref{iadmalg} to a solution of \eqref{eqftvdcs} follows from \cite[Theorem
4]{BSHE}.


\begin{thebibliography}{10}


\bibitem{BB} {\sc J. Barzilai and J.M. Borwein}, {\em
Two point step size gradient method}, IMA J. Numer. Anal., 8 (1988),
141-148.

\bibitem{AB}  {\sc A. Beck and M. Teboulle}, {\em A fast iterative shrinkage-thresholding
algorithm for linear inverse problems}, SIAM J. Imaging Sci.,
2(2009), 183-202.

\bibitem{TWIST}  {\sc J. Bioucas-Dias, and M. Figueiredo}, {\em A new
TwIst: Two-step iterative thresholding algorithm for image
restoration}, IEEE Trans. Image. Precess., 16 (2007), 2992-3004.

\bibitem{CAI1}  {\sc J.F. Cai, S. Osher, and Z. Shen}, {\em Linear
Bregman iterations for compressed sensing}, Math. Comput., 78(2009),
1515--1536.
\bibitem{CAI2}  {\sc J.F. Cai, S. Osher, and Z. Shen}, {\em Convergence of the linearized
Bregman iteration for $\ell_1$-norm minimization}, Math. Comput.,
78(2009), 2127--2136.

\bibitem{CANDES1}  {\sc E. Cand\`{e}s, J. Romberg, and T. Tao}, {\em Stable signal recovery
from imcomplete and inaccurate information}, Communications on Pure
and Applied Math., 59 (2005), 1207-1233.

\bibitem{CANDES2}  {\sc E. Cand\`{e}s, J. Romberg, and T. Tao}, {\em
Robust uncertainty principles: Exact signal reconstruction from
highly incomplete frequence information}, IEEE Trans. Inform.
Theory, 52 (2006), 489-509.

\bibitem{PLC}  {\sc P.L. Combettes and J.C. Pesquet}, {\em Proximal
thresholding algorithm for minimization over orthonormal bases},
SIAM J. Optim., 18 (2008), 1351-1376.


\bibitem{TFCHAN1}  {\sc T.F. Chan and K. Chen}, {\em An
optimization-based multilevel algorithm for tatal variation image
denoising}, Multiscal Model. Simul., 5 (2006), 615-645.

\bibitem{TFCHAN2}  {\sc T.F. Chan, G.H. Golub, and P. Mulet}, {\em A
nonlinear diffusivity fixed point method in total variation based
image restoration}, SIAM J. Sci. Comput., 20 (1999), 1964-1977.

\bibitem{DONOBO1}  {\sc D. Donoho}, {\em Compressed sensing}, IEEE Trans. Inform.
Theory, 52 (2006), 1289-1306.

\bibitem{SHARE}  {\sc D. Donoho}, {\em For most large underdetemind
systems of linear equations, the minimal l1-norm solution is also
the sparsest solution}, Communications on Pure and Applied
Mathematics, 59 (2006), 907-934.

\bibitem{ME06} {\sc M. Elad}, {\em Why simple shrinkage is still relevant for redundant representations?},
IEEE Transactions on Information Theory, 52 (2006), pp. 5559--5569.

\bibitem{MEBM06} {\sc M. Elad, B. Matalon, and M. Zibulevsky}, {\em Image denoising with shrinkage
and redundant representations}, in Proc. IEEE Computer Society
Conference on Computer Vision and Pattern Recognition, New York,
2006.

\bibitem{FNW1}  {\sc M. Figueriedo, R. Nowak, and S.J. Wright}, {\em Gradient projection for sparse
reconstruction, Application to compressed sensing and other inverse
problems}, IEEE Journal of Selected Topics in Signal Processing:
Special Issue on Convex Optimization Methods for Signal Processing,
1 (2007), 586-598.

\bibitem{GABAYM}  {\sc D. Gabay and B. Mercier}, {\em A dual
algorithm for the solution of nonliear variational problems via
finite element approximation}, Computer Math. Appl., 2 (1976),
17-40.

\bibitem{GLOWINSKI1}  {\sc R. Glowinski and P. Le Tallec}, {\em
Augmented Lagrangian and operator-splitting methods}, in: Nonlinear
Mechanice, SIAM Studies in Applied Mathematics, Philadephia, PA,
1989.

\bibitem{GLOWINSKI2}  {\sc R. Glowinski}, {\em Numerical methods for
nonlinear variational problems}, Springer-Verlat, New York, 1984.

\bibitem{ZHANG2}  {\sc E.T. Hale, W. Yin, and Y. Zhang}, {\em A fixed-point
continuation method for l1-regularized minimization with
applications to compressed sensing}, SIAM J. Optim., 19 (2008),
1107-1130.

\bibitem{BSHE}  {\sc B. He, L.Z. Liao, D. Han, and H. Yang}, {\em A
new inexact alteratin directions method for monotone variational
inequalities}, Math. Program., 92 (2002), 103-118.

\bibitem{NG3}  {\sc Y. Huang, M.K. Ng, and Y.W. Wen}, {\em  A fast total
variation minimization method for image restoration}, Multiscale
Model. Simul., 7 (2008), 775-795.

\bibitem{MRI}  {\sc M. Lustig, D. Donoho, and J. M. Pauly}, {\em Sparse MRI: the
application of compressed sensing for rapid MR Imaging}, Magnetic
Resonance in Medicine, 58 (2007), 1182-1195.

\bibitem{NG1}  {\sc M. Ng, L. Qi, Y. Yang, and Y. Huang}, {\em On
semismooth Newton's methods for total variation minimization}, J.
Math. Imaging Vision, 27 (2007), 265-276.

\bibitem{OSHER2} {\sc  S. Osher, M Burger, D. Goldfarb, J. Xu, and W.
Yin}, {\em An iterated regularization method fro total
variation-based image restoration}, Multiscale Model. Simul., 4
(2005), 460-489.

\bibitem{WYIN3}  {\sc S. Osher, Y. Mao, B. Dong, and W. Yin}, {\em
Fast linearized Bregman iteration for compressed sensing and sparse
denoising}, TR08-07, CAAM, Rice University.

\bibitem{PCF}  {\sc N. Paragios, C. Chen, and O. Faugeras}, {\em
Handbook of mathematical models in computer vision}, Springer, New
York, 2006.

\bibitem{TV}  {\sc L. Rudin, S. Osher, and E. Fatemi}, {\em Nonlinear total variation based
noise removal algorithms}, Phys. D, 60 (1992), 259-268.

\bibitem{JLMN03} {\sc J. L. Starck, M. Nguyen, and F. Murtagh}, {\em Wavelets and curvelets for
image deconvolution: a combined approach}, Signal Processing, 83
(2003), pp. 2279--2283.

\bibitem{TIKHONOV}  {\sc A. Tikhonov and V. Arsenin}, {\em
Solution of Ill-Posed Problems}, Winston, Washington, DC, 1977.

\bibitem{SPGL1}  {\sc E. van den Berg and M.P. Friedlander}, {\em Probing the
pareto frontier for basis pursuit solutions}, SIAM J. Sci. Comput.,
31 (2008), 890-912.

\bibitem{VCR96} {\sc C. R. Vogel and M. E. Oman}, {\em Iterative methods for
total variation denoising}, SIAM J. Sci. Comput., 17 (1996),
pp.227--238.

\bibitem{AF98} ---, {\em A fast, robust total variation based
reconstruction of noisy, blurred images}, IEEE Trans. Image
Process., 7 (1998), pp. 813--824.

\bibitem{YANG1}  {\sc Y. Wang, J. Yang, W. Yin, and Y. Zhang}, {\em A new alternating
minimization algorithm for total variation image reconstruction},
SIAM J. Imaging Sci., 1 (2008), 248-272.

\bibitem{NG4}  {\sc Y.W. Wen, M.K. Ng, and W.K. Ching}, {\em Iterative algorithms
based on decoupling of deblurring and denoising for image
restoration}, SIAM J. Sci. Comput., 30 (2008), 2655-2674.

\bibitem{WEN}  {\sc Z. Wen, W. Yin, D. Goldfard, and Y. Zhang}, {\em
A fast algorithm for sparse rescontruction based on shringkage:
subspace optimization and continuations}, TR09-01, CAAM, Rice
University.

\bibitem{YANG2}  {\sc J. Yang, W. Yin, Y. Zhang, and Y. Wang},
{\em A fast algorithm for edge-preserving variational multichannel
image restoration}, SIAM J. Imaging Sci., 2 (2009), 569-592.

\bibitem{YANG5}  {\sc J. Yang and Y. Zhang}, {\em Alternating direction
algorithms for $\ell_1$-problems in compressive sensing}, TR09-37,
CAAM, Rice University.

\bibitem{YANG3}  {\sc J. Yang, Y. Zhang, and W. Yin}, {\em An efficient
TVL1 algorithm for deblurring multichannel images corrupted by
impulsive noise}, SIAM J. Sci. Comput., 31 (2009), 2842-2865.

\bibitem{YANG4}  {\sc J. Yang, Y. Zhang, and W. Yin}, {\em A fast
TVL1-L2 minimization algorithm for signal reconstruction from
partial Fourier data}, IEEE J. Special Topics Signal Processing, to
appear.

\bibitem{WYIN2}  {\sc W. Yin, S. Osher, D. Goldfard, and J. Darbon},
{\em Bregman Iterative Algorithms for l1-Minimization with
Applications to compressed sensing}, SIAM J. Imaging Sci., 1 (2008),
143-168.

\bibitem{YZHANG1}  {\sc Y. Zhang}, {\em On the theory of compressed
sensing via $\ell_1$-minimization: simple derivations and
extensions}, TR08-11, CAAM, Rice University.


\end{thebibliography}
\end{document}